\documentclass[12pt]{amsart}
\usepackage{mathrsfs}
\usepackage{latexsym,amssymb,amsmath,amsfonts,graphicx,epstopdf,subfigure,float}
\usepackage{pifont}
\usepackage{cite}
\usepackage{appendix}

\def\dist{\mathop{\rm dist}\nolimits}

\newtheorem{thm}{Theorem}[section]
\newtheorem{lem}{Lemma}[section]
\newtheorem{prop}{Proposition}[section]
\newtheorem{coro}{Corollary}[section]

\newtheorem{defn}{Definition}[section]
\newtheorem{exam}{Example}[section]

\numberwithin{equation}{section}
\newtheorem{remark}{Remark}[section]

\newcommand{\bx}{{\mathbf{x}}}
\newcommand{\by}{{\mathbf{y}}}
\newcommand{\bz}{{\mathbf{z}}}
\newcommand{\ba}{{\mathbf{a}}}
\newcommand{\bb}{{\mathbf{b}}}

\newcommand{\bu}{{\mathbf{u}}}
\newcommand{\bv}{{\mathbf{v}}}

 \usepackage{color}

\begin{document}
\title{Finite state automata and homeomorphism of self-similar sets}

\author{Liangyi Huang}
\address{College of Mathematics and Statistics, Chongqing University, Chongqing, 401331, China}
\email{liangyihuang@cqu.edu.cn}

\author{Zhiying Wen} \address{Department of Mathematical Sciences, Tsinghua University, Beijing, 100084, China}
\email{wenzy@tsinghua.edu.cn}

\author{Yamin Yang$^*$}
\address{Institute of applied mathematics, College of Science, Huazhong Agricultural University, Wuhan,430070, China}
\email{yangym09@mail.hzau.edu.cn}

\author{Yunjie Zhu} \address{School of Mathematics and Physics, Hubei Polytechnic University, Huangshi, 435003, China}
\email{yjzhu\_ccnu@sina.com}

\thanks{* The correspondence author.}

\maketitle

\begin{abstract}
The topological and metrical  equivalence of fractals is an important topic in analysis. In this paper,
we use a class of finite state automata, called $\Sigma$-automaton, to construct psuedo-metric spaces, and then apply them to the study of classification of self-similar sets.
We first introduce a notion of topology automaton of a fractal gasket, which is a simplified version
of neighbor automaton; we show that a fractal gasket is homeomorphic to  the  psuedo-metric space induced by  the topology automaton. Then we construct a universal map to show that psuedo-metric spaces induced by different automata
can be bi-Lipschitz equivalent. As an application, we obtain a rather general sufficient condition for two fractal gaskets to be homeomorphic or Lipschitz equivalent.
\end{abstract}

\section{ Introduction}

To determine whether two  fractal sets are homeomorphic, quasi-symmetric or Lipschitz equivalent is important
 in analysis.
The study of homeomorphism of fractal sets dated back to
 Whyburn \cite{Why58}. For  studies of  quasi-symmetric equivalence of fractal sets,
  see \cite{Solomyak10, Bonk}.
The study of Lipschitz equivalence of fractal sets derives from 1990's and it becomes a very active topic in recent years \cite{DS,FM,FanRZ15,LuoL13,RRX06,RuanWX14,RaoZ15,XiXi20}, where most of the studies  focus on self-similar sets which are totally disconnected.

For self-similar sets which are not totally disconnected, to construct homeomorphisms, quasi-symmetric maps or bi-Lipschitz maps  is very difficult and there are few results
(\cite{Why58,  Bonk, YZ18, RaoZhu16}).
  Whyburn \cite{Why58} proved that all the Sierpinski curves are homeomorphic, which can be applied to a class of connected fractal squares.
  Solomyak \cite{Solomyak10} proved that a Julia set is always  quasi-symmetric equivalent to a planar self-similar set with two branches.
   Bonk and Merenkov\cite{Bonk} proved that the quasi-symmetric map  from  Sierpinski carpet to itself must be an isometry.

\begin{figure}[H]
\centering
\subfigure[$F_1$]{
\includegraphics[width=3.3 cm]{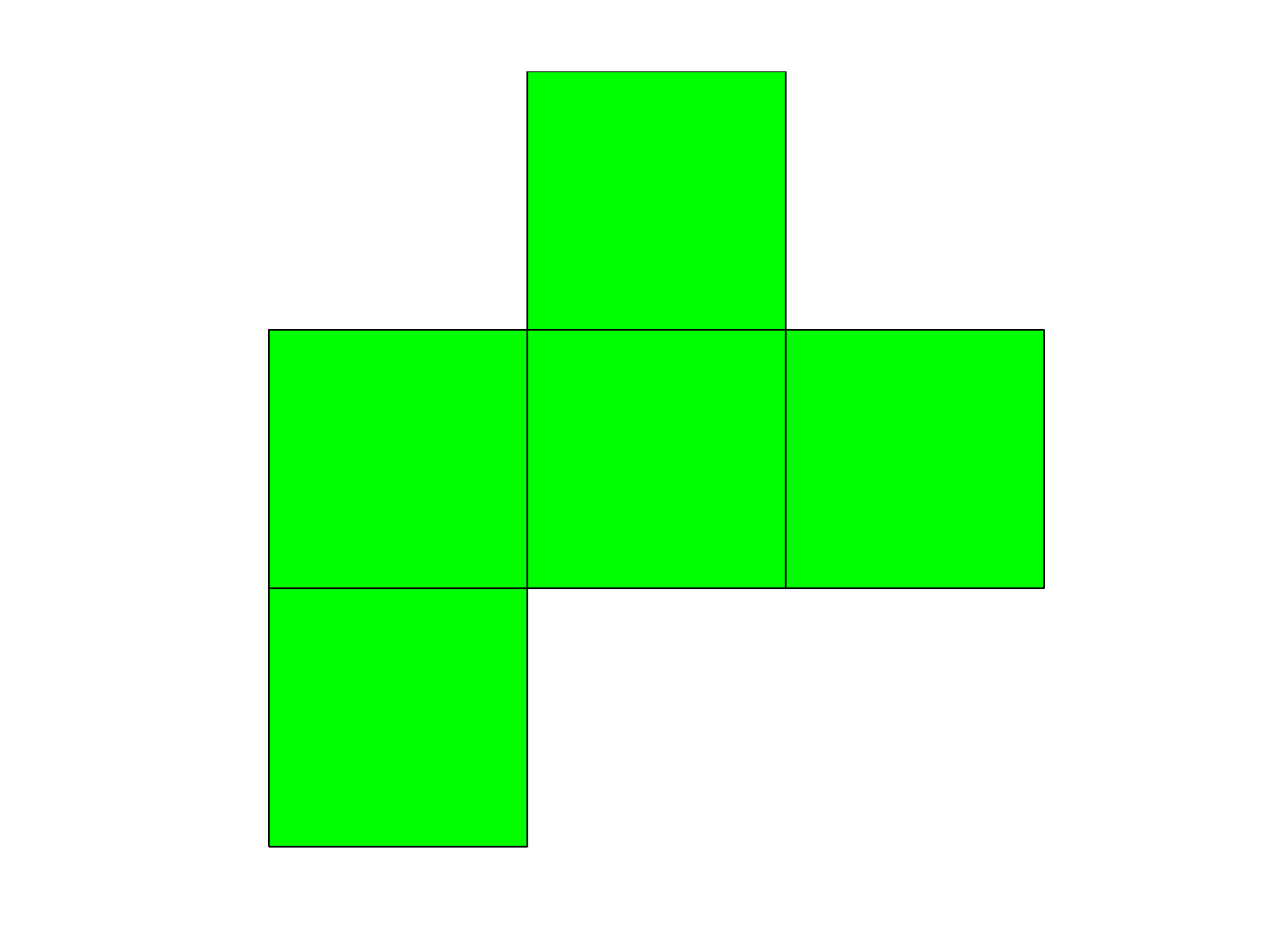}
\includegraphics[width=3.3 cm]{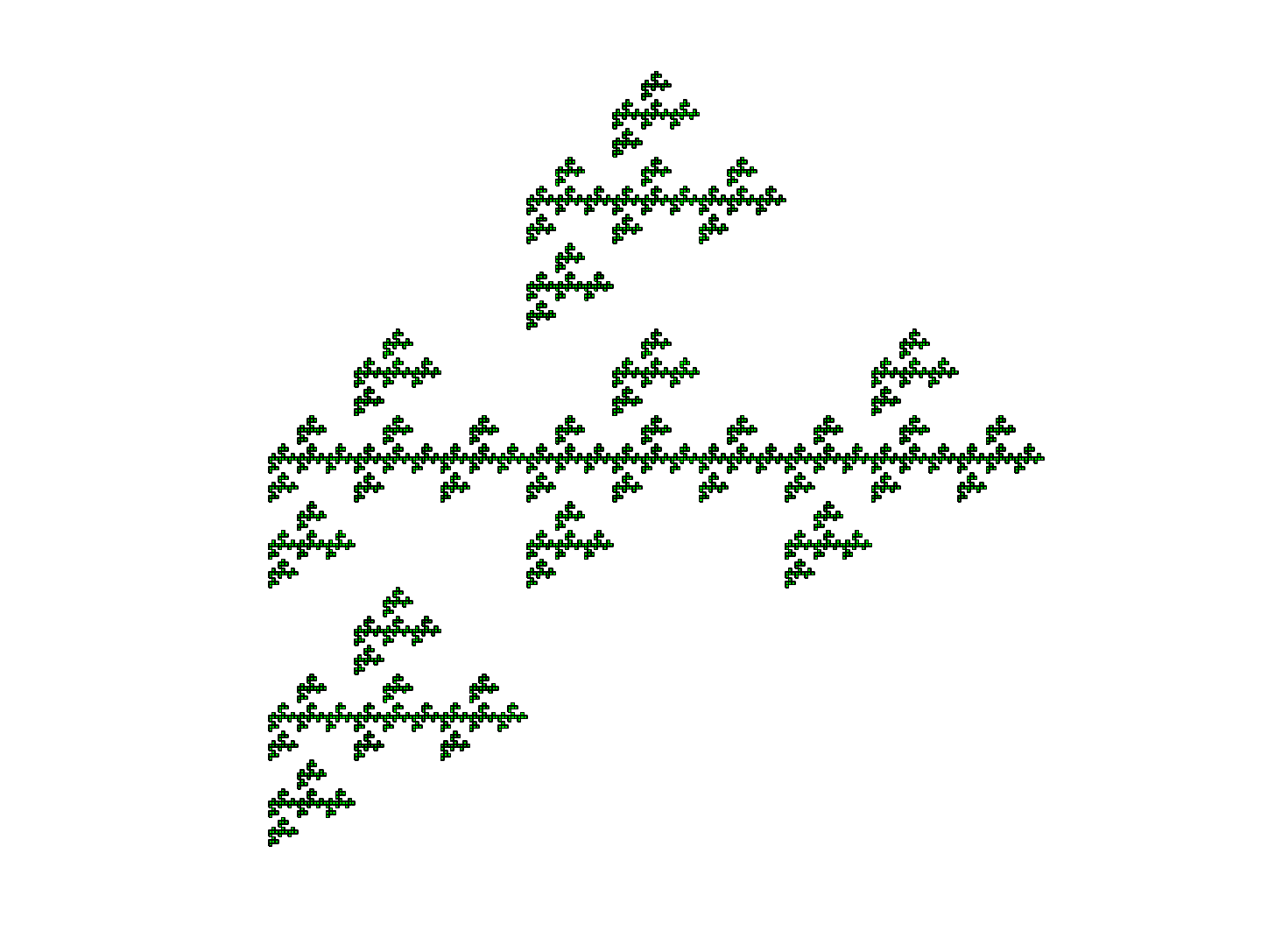}
}\quad
\subfigure[$F_2$]
{\includegraphics[width=3.3 cm]{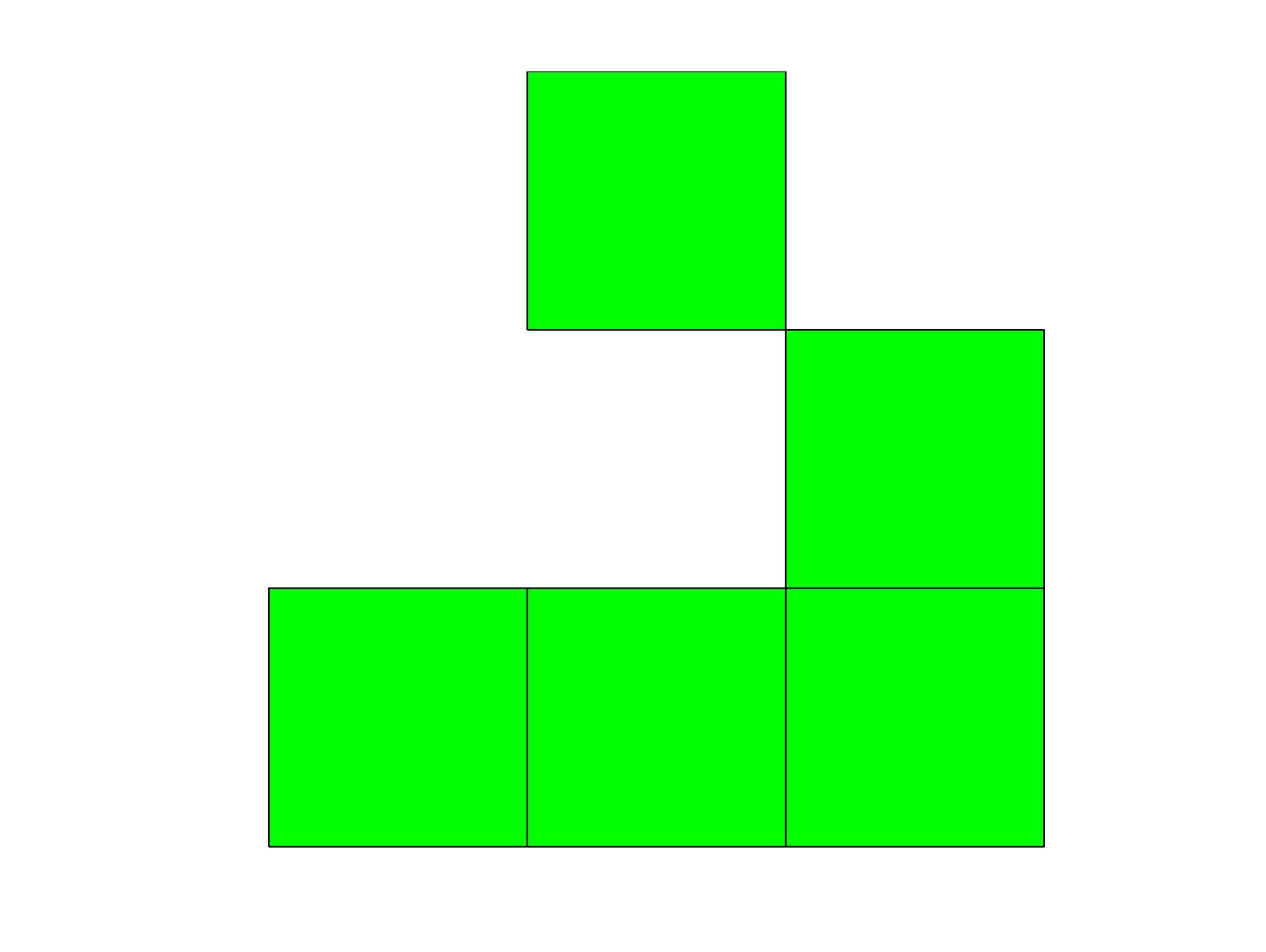}
\includegraphics[width=3.3 cm]{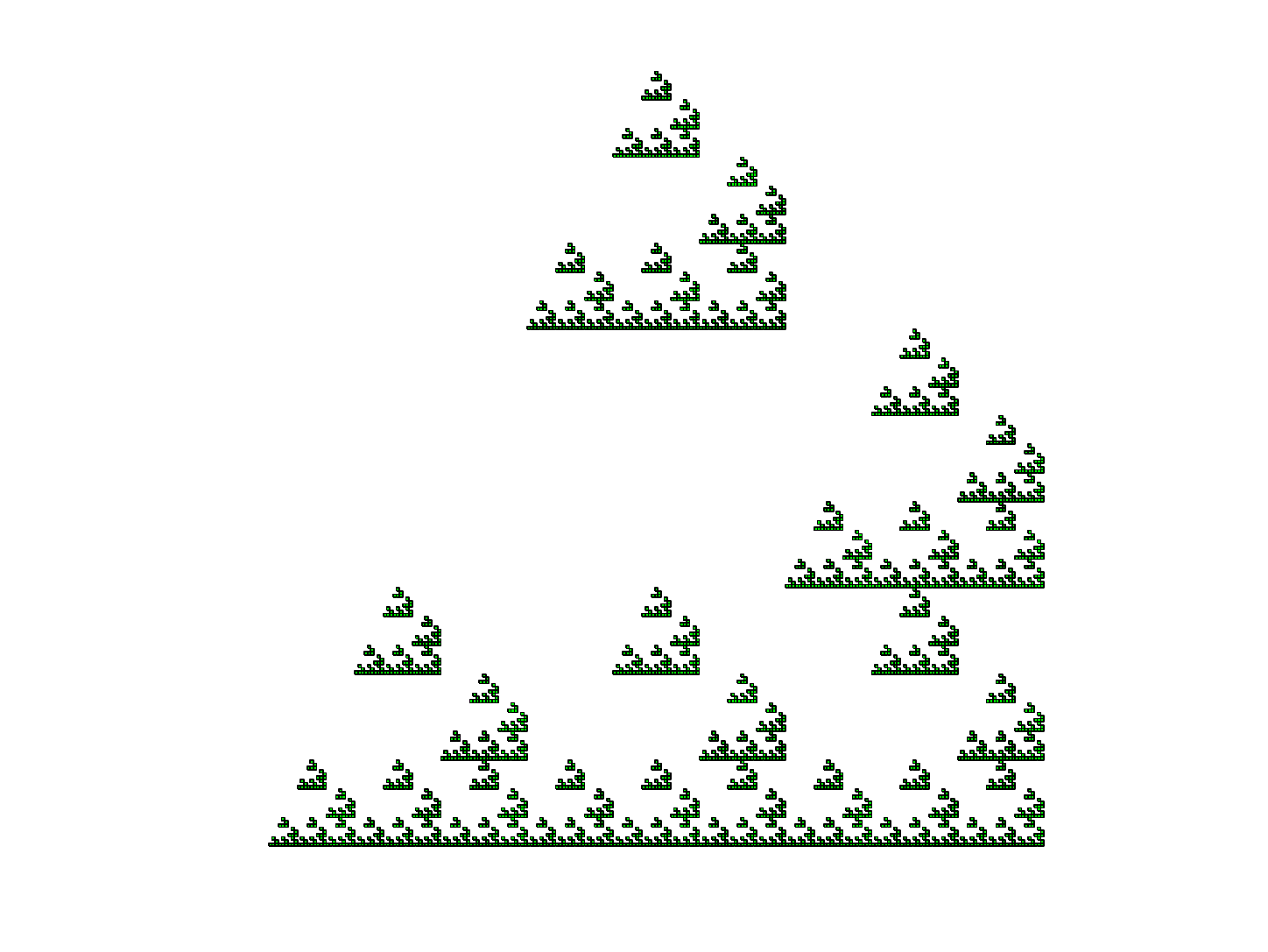}}
\\
\subfigure[$F_3$]{
\includegraphics[width=3.3 cm]{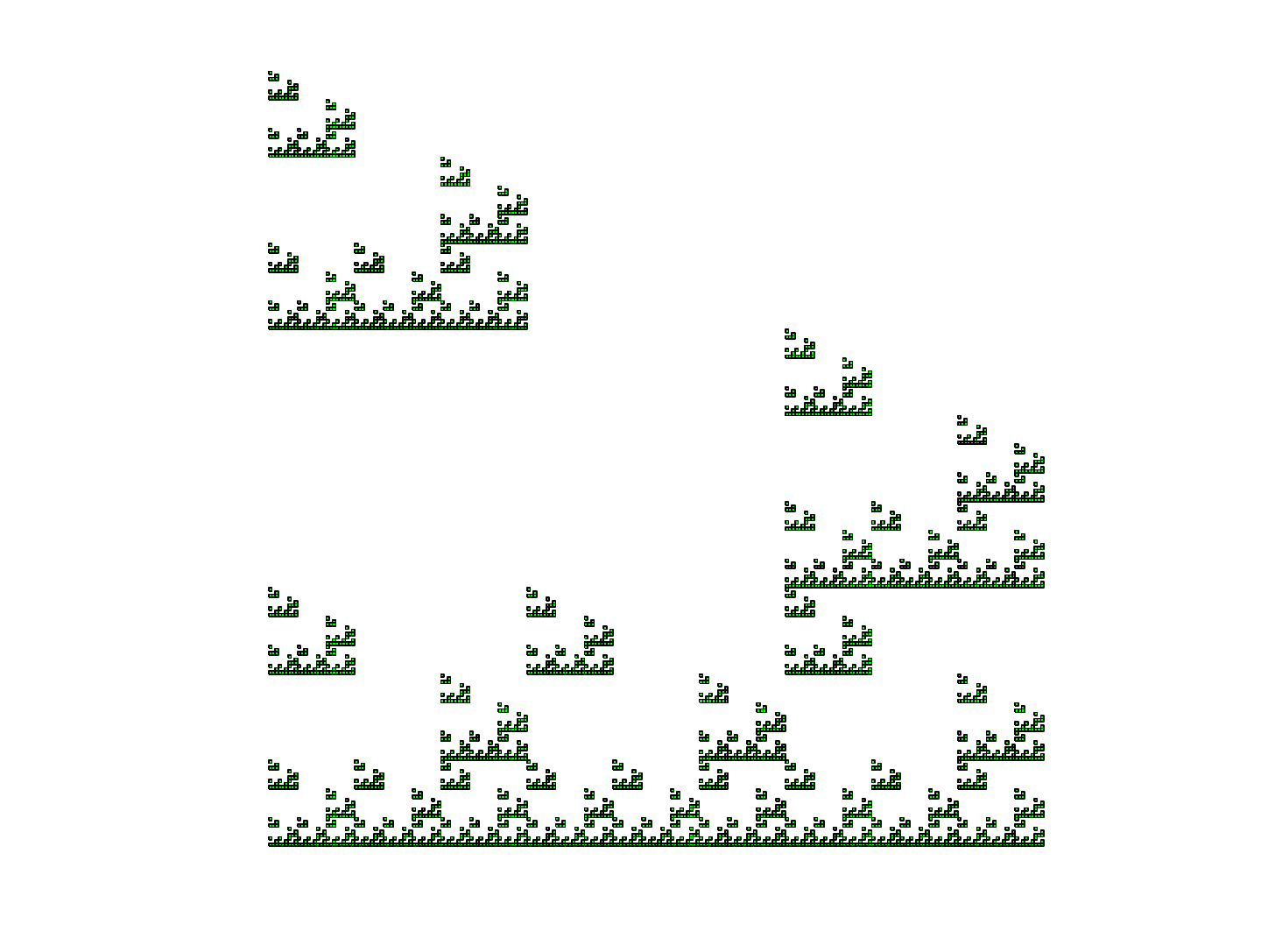}}\quad
\subfigure[$F_4$]{
\includegraphics[width=3.3 cm]{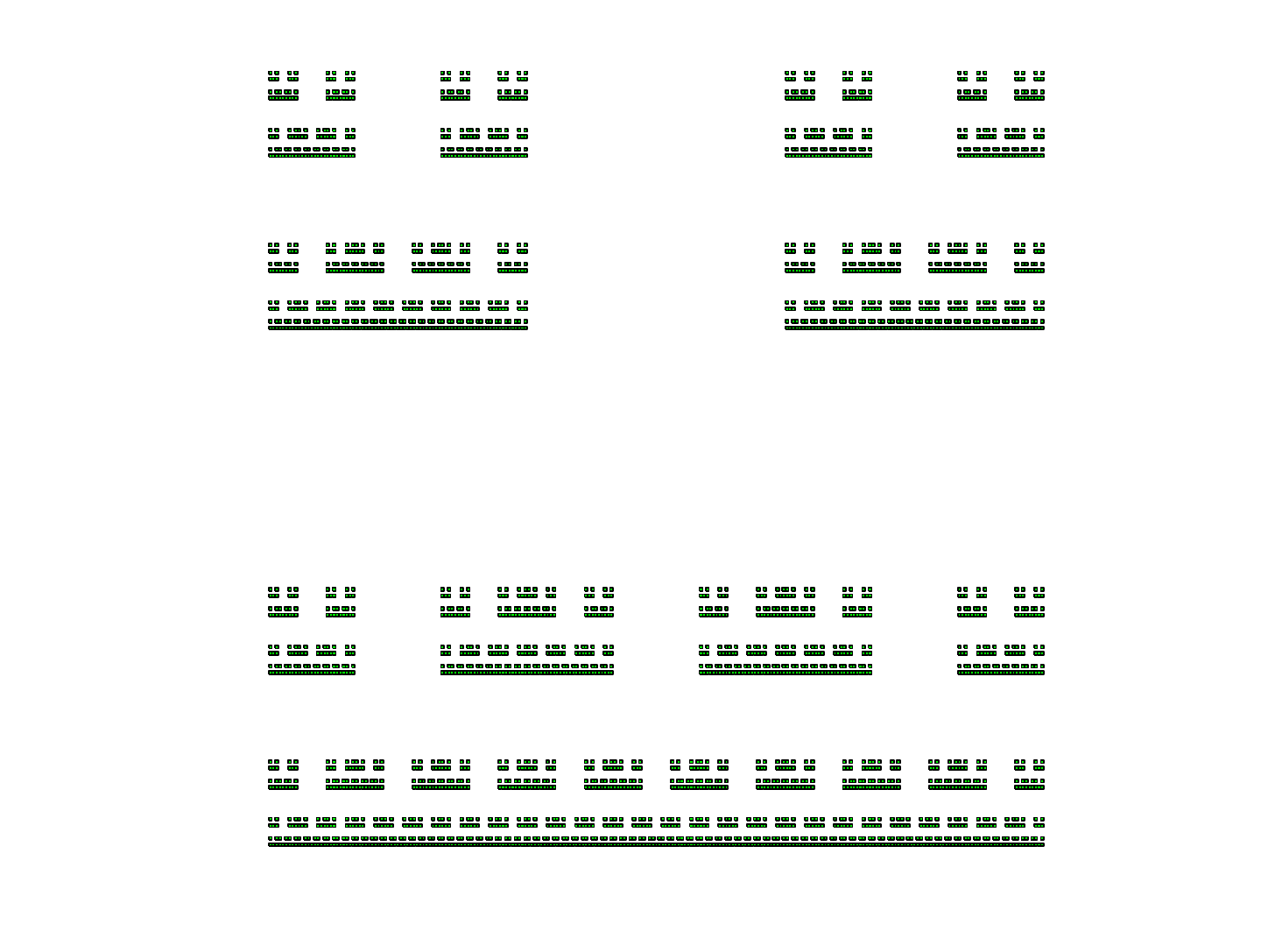}}\quad
\subfigure[$F_5$]{
\includegraphics[width=3.3 cm]{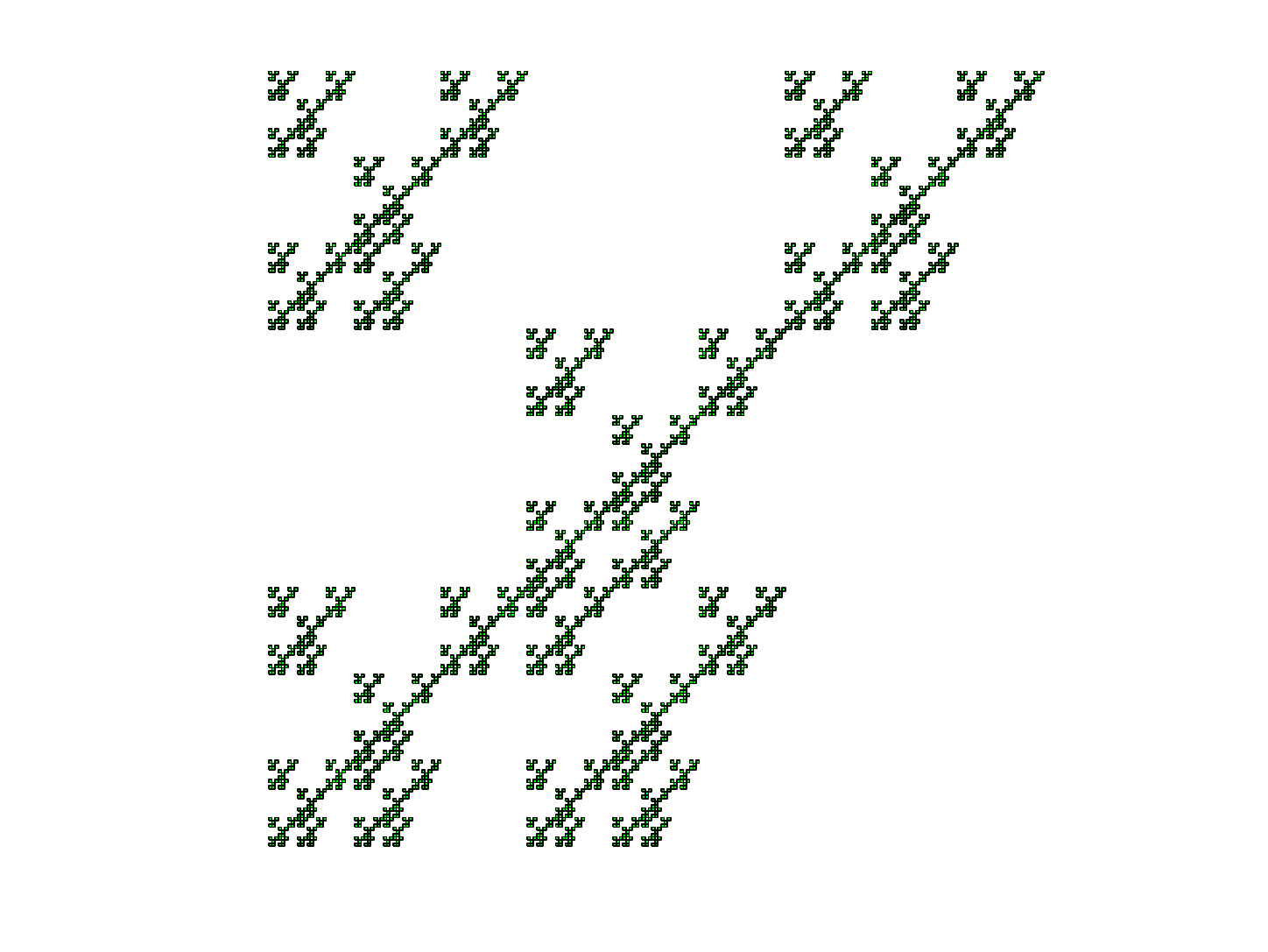}}
\\
\caption{Fractal squares with $5$ branches.}
\label{compare}
\end{figure}

There are several works devoted to
  the Lipschitz classification  of non-totally disconnected  fractal squares with contraction ratio $1/3$, that is, a kind of Sierpinski carpets (\cite{RuanW17,LuoL16,RWW17,YZ18}), but the problem is unsolved in  case
  of the fractal squares with 5 branches.
  Using  neighbor automaton,   Rao-Zhu \cite{RaoZhu16} proved that $F_1\simeq F_2$ in Figure \ref{compare}, but it is not known whether $F_j, j=2,3,4,5$  are Lipschitz equivalent or homeomorphic.

In this paper,  we develop a more
  systematic theory to  study the topological and metrical  equivalence of self-similar sets
  by finite state automaton, which are mostly inspired by Rao and Zhu \cite{RaoZhu16}.
First, we recall the definition of finite state automaton.

\begin{defn}[\cite{JEH79}]\label{FA}
\emph{A \emph{finite state automaton} is a 5-tuple $(Q,\mathcal{A},\delta,q_0,P)$, where $Q$ is a finite set of states, $\mathcal{A}$ is a finite input alphabet, $q_0$ in $Q$ is the initial state, $P\subset Q$ is the set of final states, and $\delta$ is the transition function mapping $Q\times\mathcal{A}$ to $Q$. That is, $\delta(q,a)$ is a state for each state $q$ and input symbol $a$.}
\end{defn}

 Let  $\Sigma=\{1,\dots,N\}$, which we call the set of alphabet. Inspired by the neighbor automaton of self-similar sets, we define

  \begin{defn}\emph{A finite state automaton $M$ is called a \emph{$\Sigma$-automaton} if
\begin{equation}\label{sigma-auto}
M=(Q,\Sigma^2,\delta,Id, Exit),
\end{equation}
 where the state set is $Q=Q_0\cup \{Id, Exit\}$, the input alphabet is $\Sigma^2$,  the initial state is $Id$, the final state is $Exit$, and the transition function $\delta$ satisfies
 \begin{equation}\label{eq:id}
  \delta(Id,(i,j))=Id \Leftrightarrow i=j.
 \end{equation}
 }
 \end{defn}

  By inputting symbol string $(\bx,\by)\in\Sigma^{\infty}\times\Sigma^{\infty}$ to $M$, we obtain a sequence of states $(S_{i})_{i\ge 0}$ and call it the \emph{itinerary} of $(\bx,\by)$. If we arrive at the state $Exit$, then we stop there and
 the itinerary is finite, otherwise,  it is infinite.
 We    define the \emph{surviving time} of $(\bx,\by)$ to be
\begin{equation}
T_M(\bx,\by)=(\text{length of the itinerary})-1.
\end{equation}
 Let $0<\xi<1$, we define a function $\rho_{M,\xi}$ on $\Sigma^\infty\times \Sigma^\infty$ as
\begin{equation}\label{eq:metric}
\rho_{M,\xi}(\bx, \by)=\xi^{T_{M}(\bx, \by)}.
\end{equation}

We are interested in the $\Sigma$-automaton such that $(\Sigma^\infty, \rho_{M,\xi})$ is a psuedo-quasimetric space for some $0<\xi<1$ (see Section 2 for precise definition). In this case,
we define $\bx\sim \by$ if $\rho_{M,\xi}(\bx,\by)=0$, then $\sim$ is an equivalent relation. Set
$${\mathcal A}_M=\Sigma^\infty/\sim,$$
 then $(\mathcal{A}_M,\rho_{M,\xi})$ is a psuedo-metric space (Lemma \ref{lem:induce}), and we call it the \emph{psuedo-metric space induced by $M$}.

Our purpose is to use such spaces to study the homeomorphism or Lipschitz equivalence of fractal sets.
(We remark that to construct bi-Lipschitz maps between two totally disconnected self-similar sets,  a crucial idea is to
 use  symbolic spaces  with suitable metric as  intermediate metric spaces, see \cite{FM, RRX06, XiXi10, RuanWX14, RaoZ15};  in  Luo and Lau \cite{LuoL13}, even    a certain
hyperbolic tree is used  for such purpose.)

\begin{defn}\label{def:Holder}\emph{Two psudo-metric spaces $(X,d_X)$ and $(Y,d_Y)$ are said to be \emph{H\"{o}lder equivalent},   if there is a bijection $f:~X\to Y$, a number $s>0$  and  a constant $C>0$ such that
\begin{equation}
C^{-1}d_X(x_1,x_2)^{1/s} \leq d_Y \big( f(x_1),f(x_2)\big ) \leq C d_X(x_1,x_2)^s,\ \forall x_1,x_2\in X;
\end{equation}
in this case we say $f$ is a \emph{bi-H\"{o}lder map} with index $s$. }

\emph{If $s=1$, we say $X$ and $Y$ are \emph{Lipschitz equvalent}, denote by $X\simeq Y$,  and call $f$ a \emph{bi-Lipschitz map}.
}
\end{defn}

 In this paper, we confine ourselves to a special class of  finite state automaton  called the \emph{gasket automaton},
 where the state set $Q$  contains eight states,  see  Definition \ref{triangle-auto} and \ref{gasket-auto}.
We show that

 \begin{thm}\label{lem:feasible}
 Let $M$ be a gasket automaton. Then $(\Sigma^\infty, \rho_{M,\xi})$ is a psuedo-metric space for every $0<\xi<1$.
 \end{thm}

To construct bi-H\"older maps between induced psuedo-metric spaces of different automata is a difficult problem. To this end, we define a \emph{$\gamma$-isolated} condition.
 Moreover, for a  gasket automaton $M$, we define the one-step  simplification of $M$
(see Section \ref{sec:gamma}).
The key result of this paper is the following.

\begin{thm}\label{spaceLip} Let $M$ be a  gasket automaton satisfying the $\gamma$-isolated condition, and let $M'$ be a one-step simplification of $M$.
Then $(\mathcal{A}_M,\rho_{M,\xi})\simeq(\mathcal{A}_{M'},\rho_{M',\xi})$ for every $\xi\in(0,1)$, \textit{i.e},
they are Lipschitz equivalent.
\end{thm}

Next, we apply the above results to study the classification of  fractal gaskets defining as follows.
An \emph{iterated function system} (IFS) is a family of contractions $\{\varphi_j\}_{j=1}^N$  on $\mathbb{R}^{d}$, and the \emph{attractor} of  the IFS is the unique nonempty compact set $K$ satisfying
$K=\bigcup_{j=1}^N\varphi_j(K)$ and it is called a \emph{self-similar set} \cite{Hutchinson1981} if all $\varphi_j$ are similitudes.

Let $\triangle\subset\mathbb{R}^2$ be the regular triangle with vertexes $(0,0),$ $(1,0)$, $\omega=(1/2,\sqrt{3}/2)$.
\begin{defn}[Fractal gasket]\label{fractalgasket}
\emph{ Let $(r_1,\dots, r_N)\in (0,1)^N$ and $\{d_1,\dots, d_N\}\subset {\mathbb R}^2$.
 Let $K$ be a self-similar set  generated by the IFS
$\{\varphi_j\}_{j=1}^N$ where $\varphi_j(z)=r_j(z+d_j)$.
We call $K$ a \emph{fractal gasket} if  }

\emph{ (i) $\bigcup_{i=1}^N\varphi_i(\triangle)\subset\triangle$;}

\emph{(ii) for any $i\neq j$,    $\varphi_i(\triangle)$ and $\varphi_j(\triangle)$ can only intersect at their vertices.}
\end{defn}

\textbf{Notations of $\alpha, \beta, \gamma$}. If $(0,0)\not\in K$, we set $\alpha=-1$, otherwise, we set $\alpha$ to be the symbol in $\Sigma=\{1,\dots, N\}$ such that $\varphi_\alpha((0,0))=(0,0)$.
Similarly, we set $\beta=-2$ or  $\varphi_\beta$ is the map with fixed point $(1,0)$, and set
$\gamma=-3$  or  $\varphi_\gamma$ is the map with fixed point $\omega$. Hereafter, we will denote $(0,0), (1,0)$ and $\omega$ by
$\omega_\alpha$, $\omega_\beta$ and $\omega_\gamma$, respectively.

For a fractal gasket $K$, we introduce a notion of \emph{topology automaton} of $K$, which  we denote  by $M_K$, in Section \ref{sec:structure}.
  Comparing to the neighbor automaton of self-similar sets, the topology automaton  records  less information
   since the information of size is ignored.

Denote $r_*=\min\{r_1,\dots, r_N\}$ and $r^*=\max\{r_1,\dots, r_N\}$.
Let $\pi: \Sigma^\infty\to K$ be the well-known  coding map,
 then $\pi$ induces a bijection from $\mathcal{A}_{M_K}$ to $K$ which we still denote by $\pi$
  (see Section \ref{sec:structure} for details). Actually, $\pi$ is a bi-H\"older map.

\begin{thm}\label{thm:Holder} Let $K$ be a fractal gasket. Let $s=\sqrt{\log r^*/\log r_*}$ and $\xi=(r_*)^s$. Then $\pi: (\mathcal{A}_{M_K},\rho_{M_K,\xi})\to K$ is a bi-H\"{o}lder map with index $s$.
\end{thm}

\begin{remark}\emph{
We remark that if two of $\{\omega_\alpha,\omega_\beta,\omega_\gamma\}$ do not belong to $K$, then   $\varphi_i(K)\cap \varphi_j(K)=\emptyset$ whenever $i\neq j$.
In this case $K$ is totally disconnected and we will deal with it in another paper. From now on,
we will always assume that $\omega_\alpha,\omega_\beta\in K$ without loss of generality,
or equivalently, $\alpha, \beta\in \Sigma$.
}
\end{remark}

A fractal gasket $K$ is said to  satisfy the \emph{top isolated condition},
if  $\omega_\gamma\in K$, and  $\varphi_{\gamma}(\triangle)\cap  \varphi_j(\triangle)=\emptyset$ provided $j\neq \gamma$. This condition corresponds to the $\gamma$-isolated condition in Theorem \ref{spaceLip}.

 We remark that if a fractal gasket satisfies the top isolated condition, then
its non-trivial connected components are horizontal line segments. See Appendix A.

\begin{defn}[Horizontal block]\label{Hblock}
\emph{Let $K$ be a fractal gasket.
 We call
$$
I=\{i_1,i_2,\dots,i_k\}\subset\Sigma
$$
a \emph{horizontal block} of $K$ if $\varphi_{i_{j}}(\omega_\beta)=\varphi_{i_{j+1}}(\omega_\alpha)$
for $1\le j\le k-1$ and $I$ is maximal with this property. We call $k$ the size of $I$. }
\end{defn}

If $\alpha$ and $\beta$ belong to a same horizontal block, then we call this block the $\alpha\beta$-block.

If $K$ satisfies the top isolated condition or $\omega_\gamma\not\in K$, and $\alpha$ and $\beta$ does not belong the same horizontal block, then $K$ is totally disconnected, and it is out of our consideration.
Let ${\mathcal F}_{T,\alpha\beta}$ denote the collection of fractal gasket $K$  satisfying

(i) $K$ satisfies the top isolated condition or  $\omega_\gamma\not\in K$;

(ii)   $\alpha$ and $\beta$ belong to a same horizontal block of $K$.

\begin{thm}\label{thm:Lip}
Let $E, F\in {\mathcal F}_{T,\alpha\beta}$. If there is a size-preserving bijection from the collection of horizontal-blocks of $E$ to that of $F$,
and the $\alpha\beta$-block of $E$ have equal size with that of $F$, then $E$ is bi-H\"{o}lder equivalent (and  homeomorphic) to $F$.

If in addition that both $E$ and $F$ are of uniform contraction ratio $r$, then $E\simeq F$.
\end{thm}

\begin{figure}
  \centering
  \subfigure[$E$]{\includegraphics[width=0.34\textwidth]{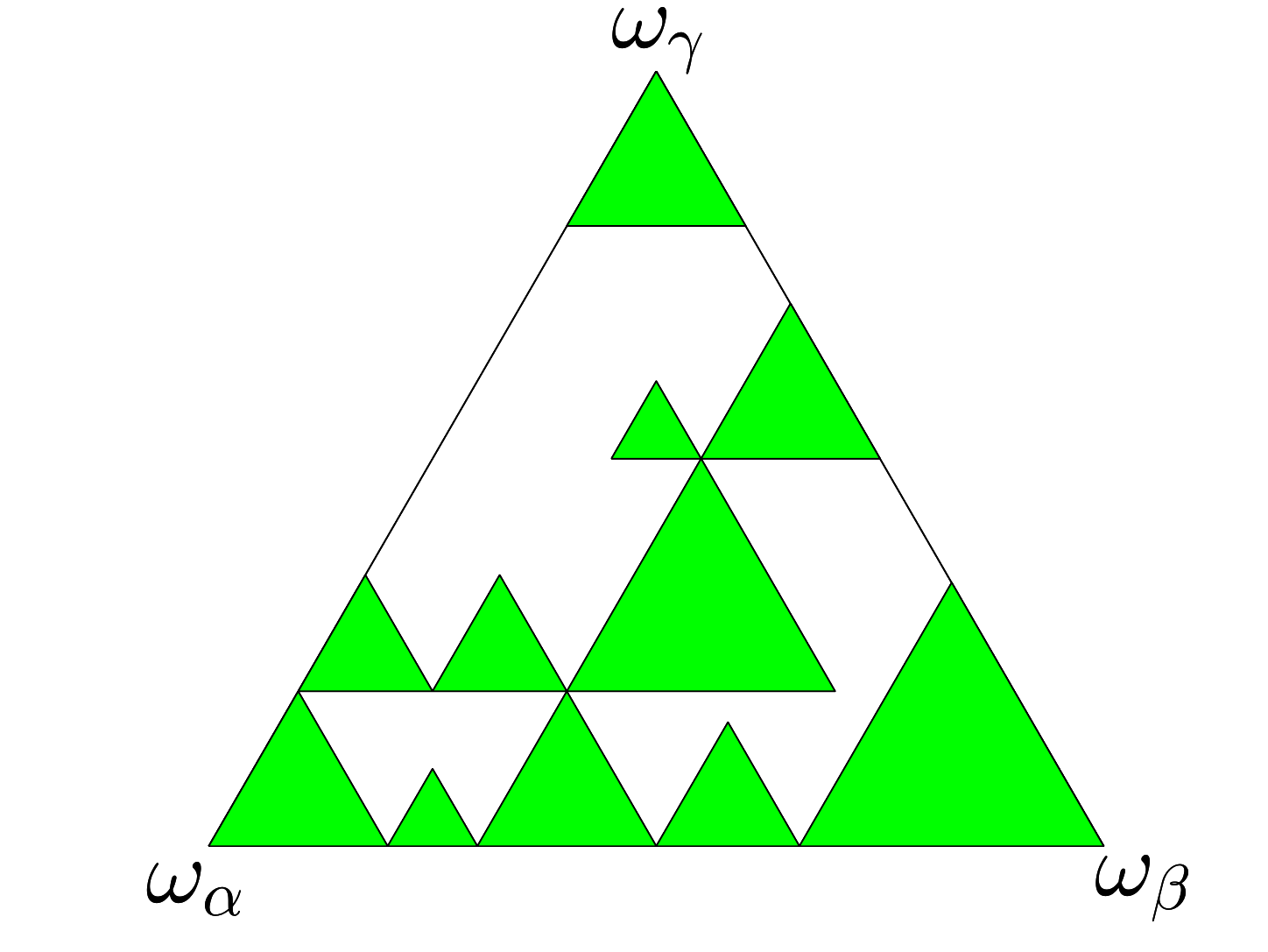}}\quad
  \subfigure[$F$]{\includegraphics[width=0.34\textwidth]{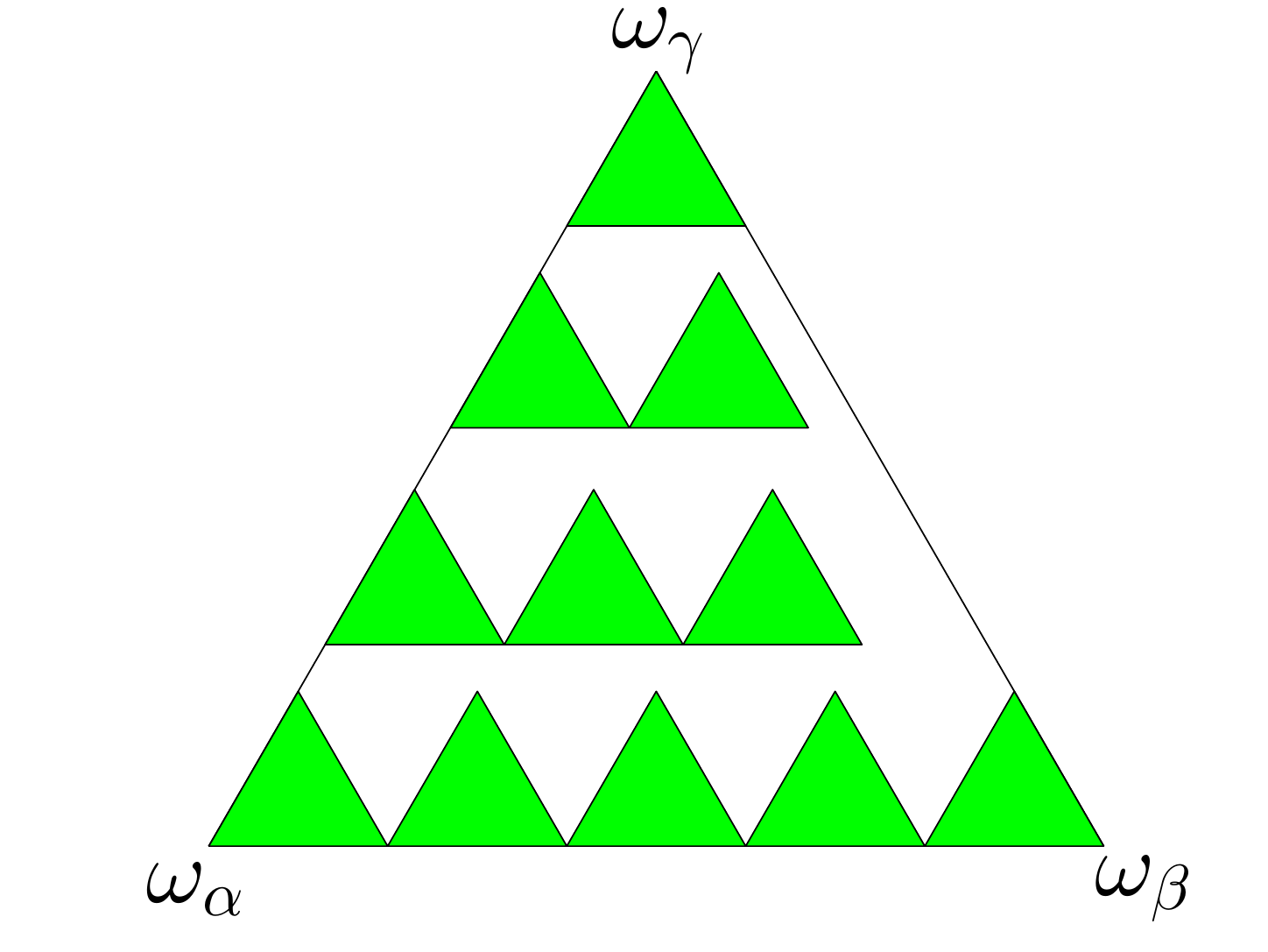}}
\caption{Fractal gaskets in Example 1.1. $E$ is homeomorphic to $F$.}
\label{example}
\end{figure}

Let us briefly explain the strategy of the proof of Theorem \ref{thm:Lip}.
Let $E$ and $F$ be two fractal gaskets in Theorem \ref{thm:Lip}. Let $M_E$  be the topology automaton of $E$.
We construct a sequence of automata
$$
M_E=M_{E,0}, M_{E,1}, \dots, M_{E,p}=M_E^*
$$
such that $M_{E,i+1}$ is a one-step simplification of $M_{E,i}$ for each $i$,
  and $M_E^*$ only records the horizontal connective relations among $\varphi_j(\triangle), j\in\Sigma$.
  By Theorem \ref{spaceLip}, we have $\mathcal{A}_{M_E}\simeq \mathcal{A}_{M_E^*}$.
  We do the same thing for $F$ and we obtain  an automaton $M_F^*$. The assumptions in Theorem \ref{thm:Lip} guarantee that $\mathcal{A}_{M_E^*}\simeq \mathcal{A}_{M_F^*}$. Therefore, $\mathcal{A}_{M_E}\simeq \mathcal{A}_{M_F}$. Finally,  by Theorem \ref{thm:Holder}, we obtain that
  $E$ is bi-H\"{o}lder equivalent to $F$.

\begin{exam}\label{exam:1}
\emph{Let $E$ and $F$ be two fractal gaskets indicated by Figure \ref{example}.
There are four horizontal-blocks in both $E$ and $F$, and the sizes of the blocks are $1, 2, 3, 5$ respectively. By Theorem \ref{thm:Lip}, $E$ is homeomorphic to $F$. (In the fractal $E$, there are 3 points connecting $\varphi_j(\triangle)$ in different horizontal blocks, so we need  $3$ one-step simplifications to obtain $M^*_E$.)}
\end{exam}

\noindent \textbf{Open question:} \emph{Can we replace the $\gamma$-isolated condition by the following condition: $\varphi_\alpha(\triangle)$ and $\varphi_\beta(\triangle)$ belong to a same connected component of
$\bigcup_{j=1}^N \varphi_j(\triangle)$, but $\varphi_\gamma(\triangle)$ does not?}

\medskip

This article is organized as follows:
In Section \ref{sec:psuedo}, we discuss the  psuedo-metric space induced by a finite state automaton. In Section 3, we introduce the gasket automaton, and Theorem \ref{lem:feasible} is proved there. In Section \ref{sec:structure}, we define the topology automaton of a fractal gasket, and Theorem \ref{thm:Holder} is proved there. In Section \ref{sec:gamma}, we  discuss the one-step simplification of a gasket automaton.
 In Section \ref{sec:proof}, we prove Theorem \ref{spaceLip} and Theorem \ref{thm:Lip}
 by assuming Theorem \ref{main}. In Section \ref{sec:mapg}, we construct a universal map $g$ on the symbolic space $\Omega$. Finally, in Section \ref{sec:time}, we prove  Theorem \ref{main} which is technical.

\section{\textbf{Psuedo-metric space induced by $\Sigma$-automaton}}\label{sec:psuedo}

Let us recall the definition of psuedo-metric space, see for instance \cite{Pepo90, Bao18}.

\begin{defn}\label{def:psuedo}\emph{
 A \emph{pseudo-quasimetric space}
is a pair $(\mathcal{A},\rho)$ where $\mathcal{A}$ is a set and $\rho:\mathcal{A}\times\mathcal{A}\rightarrow\mathbb{R}_{\geq 0}$ satisfying for all $x,y,z\in \mathcal{A}$, it holds that }

\emph{(i) $\rho(x,x)=0$;}

\emph{(ii)   $\rho(x,y)=\rho(y,x)$;}

\emph{(iii) (psuedo-triangle inequality) $\rho(x,z)\le C(\rho(x,y)+\rho(y,z)),$
  where $C\ge 1$ is a constant independent of $x,y,z$.}

\emph{If in addition $x\neq y$ implies $\rho(x,y)>0$, then we call $({\mathcal A}, \rho)$
a \emph{psuedo-metric space}. }
\end{defn}

Let $({\mathcal A},\rho)$ be a psuedo-quasimetric space.
Define $x\sim y$ if $\rho(x,y)=0$, then clearly $\sim$ is an equivalence relation.
Denote the equivalent class containing $x$ by $[x]$ . Set
$\widetilde {\mathcal A}:={\mathcal A}/\sim$
to be the quotient space.
For $[x],[y]\in\widetilde {\mathcal A}$, define
\begin{equation}\label{metric}
 \widetilde{\rho}([x],[y])=\inf\{\rho(a,b); a\in [x], b\in [y]\}.
\end{equation}

\begin{lem}\label{lem:induce}
The quotient space  $(\widetilde {\mathcal A}, \tilde \rho)$ is a psuedo-metric space.
\end{lem}

\begin{proof} The assertions $\tilde \rho([x],[x])=0$ and $\tilde \rho([x],[y])=\tilde \rho([y],[x])$ are obvious.

If $a,a'\in [x]$ and $b,b'\in [y]$, by the psuedo-triangle inequality, one can show that
$\rho(a,b)\leq C^2\rho(a',b')$. Hence, if $a\in [x]$ and $b\in [y]$, we have
\begin{equation}\label{eq:xi2}
C^{-2} \rho(a,b)\leq  \widetilde  \rho([x],[y]).
\end{equation}
It follows that $\widetilde{\rho}([x],[y])>0$ if $[x]\neq [y]$, and
\begin{equation}\label{eq:xi3}
\widetilde{\rho}([x],[z])\leq \rho(x,z)\leq C(\rho(x,y)+\rho(y,z))\leq C^3( \widetilde{\rho}([x],[y])+\widetilde{\rho}([y],[z])).
\end{equation}
The lemma is proved.
\end{proof}

Let $(\mathcal{A},\rho)$ be a pseudo-metric space. In the same manner as the metric space, we can define
convergence of sequence, dense subset and completeness of ${\mathcal A}$.  (See \cite{Pepo90, Bao18}.)
The following lemma is obvious.

\begin{lem}\label{extendLip}
Let $(\mathcal{A},\rho)$ and $(\mathcal{A}',\rho')$ be two complete pseudo-metric spaces.
Suppose  $B\subset\mathcal{A}$ is $\rho$-dense in $\mathcal{A}$ and  $B'\subset\mathcal{A}'$ is $\rho'$-dense in $\mathcal{A}'$. If $B\simeq B'$, then $\mathcal{A}\simeq\mathcal{A}'$.
\end{lem}

Let $\Sigma=\{1,\dots, N\}$.   For $a\in\Sigma$, we use $a^k$ to denote the word consisting of $k$ numbers of $a$. Let $\Sigma^{\infty}$ and $\Sigma^{k}$ be the sets of infinite words and words of length $k$ over $\Sigma$ respectively. Let $\Sigma^*=\bigcup_{k\geq0} \Sigma^{k}$.

Let $M$ be a $\Sigma$-automaton.
If $(\Sigma^\infty, \rho_M)$ is a psuedo-quasimetric space, then we can define a psuedo-metric space
by \eqref{metric}, which we denote by $({\mathcal A}_M, \rho_M)$.
Denote by $\bx\wedge\by$ the maximal common prefix of $\bx$ and $\by$. By \eqref{eq:id}, we see that
$$T_M(\bx,\by)\geq |\bx \wedge \by|,$$
 where $|W|$ denotes the length of a word $W$.

\begin{lem}\label{complete} If $(\Sigma^\infty, \rho_M)$ is a psuedo-quasimetric space, then the
induced  pseudo-metric space $(\mathcal{A}_M,\rho_M)$ is complete.
\end{lem}
\begin{proof}We equip $\Sigma^\infty$ with the following metric: For any $\bx,\by\in\Sigma^\infty$, define $d(\bx,\by)=2^{-|\bx\wedge\by|}$. It is folklore that $(\Sigma^\infty,d)$ is a compact metric space.

For any Cauchy sequence $\{[\bx_k]\}_{k=1}^\infty$ of $\mathcal{A}_M$, let $\{\bx_{k_p}\}_{p=1}^\infty$ be a subsequence of $\{\bx_k\}_{k=1}^\infty$ which converges to $\by\in\Sigma^\infty$ in the metric $d$.
Then $\lim_{p\rightarrow\infty}|\bx_{k_p}\wedge\by|=+\infty$, so we have $[\bx_{k_p}]\to[\by]$ with respect to $\rho_M$. Since $\{[\bx_k]\}_{k=1}^\infty$ is Cauchy in $\rho_M$,  we conclude that $[\bx_k]\to[\by]$ in $\rho_M$. The lemma is proved.
\end{proof}

\begin{lem}\label{Omegadense} Suppose $(\mathcal{A}_M, \rho_M)$ is a psuedo-metric space. Let $\kappa\in \Sigma$. Then the set $\Omega=\{[\omega\kappa^{\infty}];~\omega\in\Sigma^*\}$ is $\rho_M$-dense in $\mathcal{A}_M$.
\end{lem}
\begin{proof}
Pick $[\bx]\in\mathcal{A}_M$, denote $\bx=(x_i)_{i=1}^\infty$ and let $\bx_k=x_1\dots x_k\kappa^\infty$, then $[\bx_k]\in\Omega$. Clearly, $[\bx_k]\to[\bx]$. This finishes the proof.
\end{proof}

\section{\textbf{Triangle automaton and gasket automaton}}
In this section, we   introduce the triangle automaton and gasket automaton.

\subsection{Triangle automaton}
\ \\
\indent Let $\Sigma=\{1,\cdots, N\}$.
Let $\{\alpha,\beta,\gamma\}$ be a subset of $\Sigma\cup \{-1,-2,-3\}$.
For  a pair $u,v\in \{\alpha,\beta,\gamma\}$ with $u\ne v$, we associate with it a state and denote it by $S_{uv}$. We set the state set $Q$ to be
\begin{equation}\label{states}
Q=\{S_{\alpha\gamma},S_{\beta\gamma},S_{\alpha\beta},S_{\beta\alpha},
S_{\gamma\beta},S_{\gamma\alpha}\}\cup \{Id, Exit\}.
\end{equation}

\begin{defn}[Triangle automaton]\label{triangle-auto}
\emph{A $\Sigma$-automaton $M=\{Q,\Sigma^2,\delta,Id,Exit\}$ is called a \emph{triangle automaton} if
 $Q$ is given by \eqref{states}, and  the transition function $\delta$ satisfies the following conditions: Let $u,v\in\{\alpha,\beta,\gamma\}$ and $u\ne v$.}


\emph{ (i) (Symmetry) If $\delta(Id,(i,j))=S_{uv}$, then $\delta(Id,(j,i))=S_{vu}$.}



\emph{ (ii) (Loop) $\delta(S_{uv},(i,j))=\left \{
\begin{array}{ll}
S_{uv},  &\text{ if } (i,j)=(v,u);\\
Exit, & \text{ otherwise.}
\end{array}
\right .$
}
\end{defn}

Let us call $S_{uv}$ the \emph{mirror state} of $S_{vu}$, and let the mirror state of $Id$ and $Exit$ to be themselves. Clearly, if $S$ is translated to $S'$ by $(i,j)$, then the mirror state of $S$
is translated to the mirror state of $S'$ by $(j,i)$. Therefore, we have

\begin{lem}\label{lem:symmetry} Let $(S_k)_{k\geq 1}$ and $(S'_k)_{k\geq 1}$  be the itineraries of $(\bx,\by)$  and
$(\by,\bx)$, respectively. Then $S_k'$ is the mirror state of $S_k$. Consequently, $T_M(\bx,\by)=T_M(\by,\bx)$.
\end{lem}

\begin{figure}[H]
  \centering
  \includegraphics[width=10 cm]{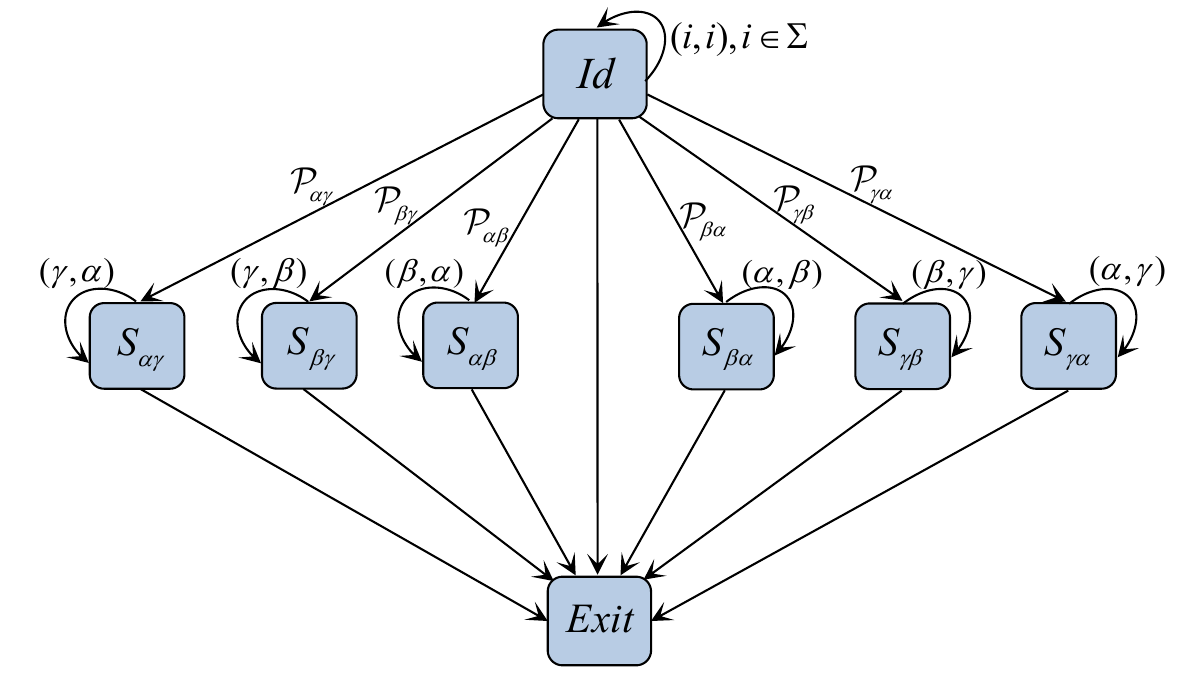}\\
\caption{The transition diagram of triangle automaton $M$.}
\label{diagram}
\end{figure}

We denote
$$
\mathcal{P}_{uv}=\{(i,j)\in\Sigma^2;\delta(Id,(i,j))=S_{uv}\}.
$$
Then by symmetry of the automaton, $\mathcal{P}_{vu}=\{(i,j);(j,i)\in\mathcal{P}_{uv}\}$.
The transition diagram of a triangle automaton $M$ is illustrated in Figure \ref{diagram}. Clearly, $M$ is completely determined by $\mathcal{P}_{\alpha\beta},\mathcal{P}_{\alpha\gamma},\mathcal{P}_{\beta\gamma}$.

If $(i,j)\in\mathcal{P}_{\alpha\gamma}$, then we denote $i\lhd_{\alpha\gamma}j$ and say $i$ is the \emph{$\alpha\gamma$-predecessor} of $j$ and $j$ is the \emph{$\alpha\gamma$-successor} of $i$; similarly, we  define $i\lhd_{\beta\gamma}j$, $i\lhd_{\alpha\beta}j$.
Moreover, according to Definition \ref{triangle-auto}(i), we make the convention that  $i\lhd_{uv}j$ if and only if $j\lhd_{vu}i$.

\subsection{Gasket automaton}
 \ \\
\indent Firstly, we recall some notions of graph theory, see \cite{Bal2000}. Let $G=(V, \mathcal{E})$ be a directed graph, where $V$ is the {vertex set} and $\mathcal{E}$ is the {edge set}. Each edge $\mathbf{e}$ is associated to an ordered pair $(u,v)$ in $V$, we say $\mathbf{e}$ is \emph{incident out} of $u$ and \emph{incident into} $v$.
 The number of edges incident out of a vertex $v$ is the \emph{outdegree} of $v$ and is denoted by $\deg^+(v)$. The number of edges incident into a vertex $v$ is the \emph{indegree} of $v$ and is denoted by $\deg^-(v)$. If $\deg^-(v)=0$, then we say $v$ is \emph{minimal}; if $\deg^+(v)=0$, then we say $v$ is \emph{maximal}. If $v$ is both minimal and maximal, then we say $v$ is \emph{isolated}.

A \emph{directed walk} joining vertex $v_1$ to vertex $v_k$ in $G$ is a sequence $(v_1, v_2, \dots, v_k)$ with $(v_{i},v_{i+1})\in\mathcal{E}$. In addition, if all $v_i (1\le i\le k)$ are distinct, then we call it a \emph{path}. If all $v_i (1\le i\le k-1)$ are distinct and $v_k=v_1$, then we call it a \emph{cycle}.
A path $(v_1, v_2, \dots, v_k)$ is called a \emph{chain},  if $v_1$ is minimal and $v_k$ is maximal.


  For  a triangle automaton $M$,  we will regard $(\Sigma,\mathcal{P}_{\alpha\gamma})$, $(\Sigma,\mathcal{P}_{\beta\gamma})$ and $(\Sigma,\mathcal{P}_{\alpha\beta})$ as three graphs.
 A symbol $j\in \Sigma$ is said to be \emph{$\alpha\beta$-minimal} (resp. maximal) if it is minimal (resp. maximal) in $(\Sigma,\mathcal{P}_{\alpha\beta})$. Similarly, we can define $\alpha\gamma$-minimal (maximal) and $\beta\gamma$-minimal (maximal).
  Moreover, we make the convention that $j$ is $uv$-minimal
 if it is $vu$-maximal.

\begin{defn}[Gasket automaton]\label{gasket-auto}
\emph{A triangle automaton $M=\{Q,\Sigma^2,\delta,Id,Exit\}$ is called a \emph{gasket automaton} if
   $\delta$ satisfies the following conditions: Let $u,v\in\{\alpha,\beta,\gamma\}$ and $u\ne v$.}

\emph{(i) (Uniqueness) If $i\lhd_{uv} j$ and $i\lhd_{uv} j'$, then $j=j'$.}

\emph{(ii) (Gathering condition) Any two of the following statements imply the third one: \ding{172} $a\lhd_{\alpha\gamma}c$; \ding{173} $a\lhd_{\beta\gamma}b$; \ding{174} $b\lhd_{\alpha\beta}c$.}

\emph{(iii)  (Boundary condition) If $\alpha\in \Sigma$, then $\alpha$ is $\alpha\gamma$-minimal and $\alpha\beta$-minimal;  similarly,  if $\beta\in \Sigma$, then  $\beta$ is $\beta\gamma$-minimal and  $\beta\alpha$-minimal; if
 $\gamma\in \Sigma$, then $\gamma$ is $\gamma\alpha$-minimal and $\gamma\beta$-minimal.}
\end{defn}

By the uniqueness property, we see that for any $u,v\in\{\alpha,\beta,\gamma\}$, the graph $(\Sigma,\mathcal{P}_{uv})$ is a union of disjoint chains and cycles.

\begin{figure}[H]
  \centering
  \includegraphics[width=6 cm]{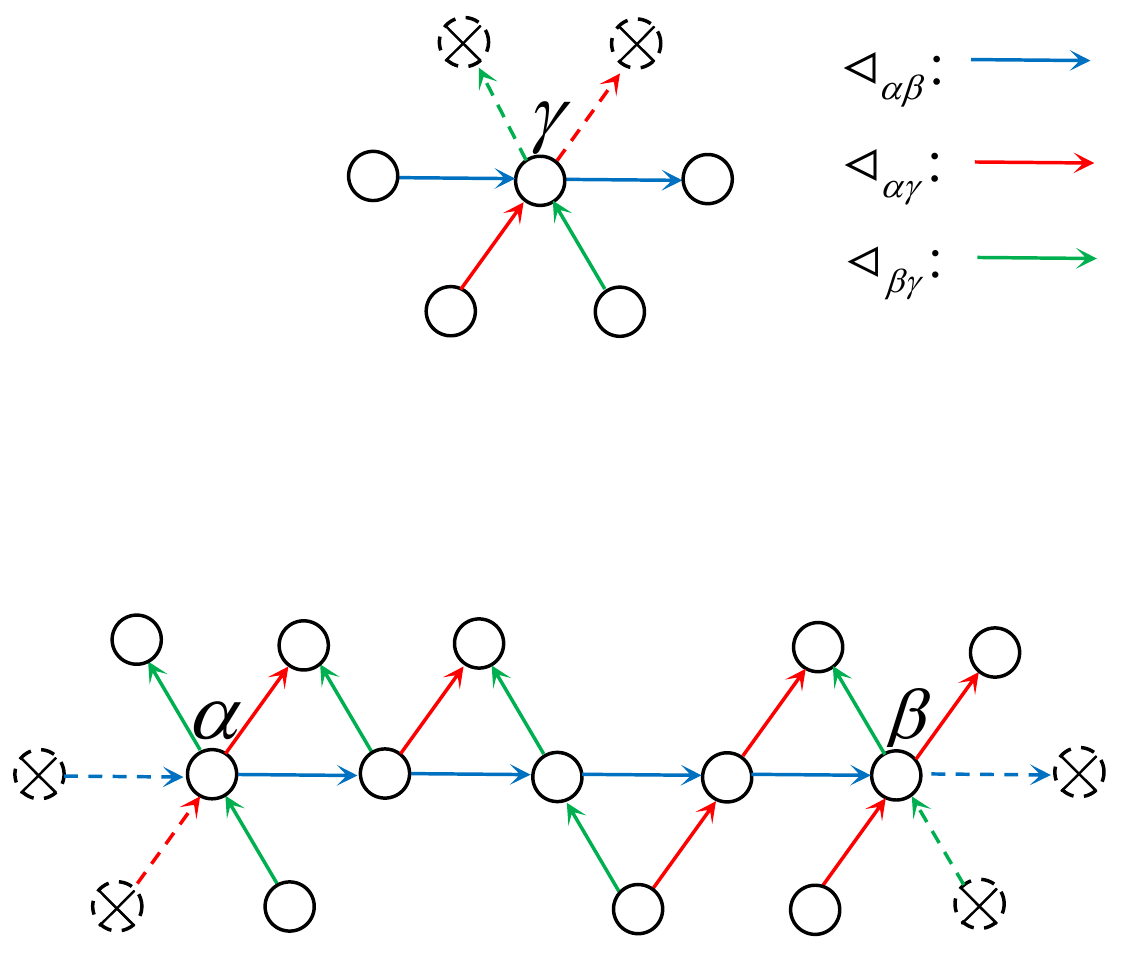}\\
\caption{An illustration  of the boundary condition: the six positions with crosses are forbidden.}
\label{fig_boundary}
\end{figure}

\begin{remark} \emph{Let $K$ be a fractal gasket. If  $\varphi_a(\omega_\gamma)=\varphi_b(\omega_\beta)=\varphi_c(\omega_\alpha)$, then
$\varphi_a(\triangle)$, $\varphi_b(\triangle)$ and $\varphi_c(\triangle)$ gather at one point.
This is the motivation of gathering condition.
}
\end{remark}

\subsection{Induced psuedo metric space}
 \ \\
\indent Let $M$ be a triangle automaton.
For $(\bx,\by)\in\Sigma^{\infty}\times\Sigma^{\infty}$,
 let $(S_{M,i})_{i\ge 0}$ be the itinerary of $(\bx,\by)$; recall that the surviving time
 $T_M(\bx,\by)$ is the largest $k$ such that $S_{M,k}\neq Exit.$

 We will use $S\stackrel{(i,j)}\longrightarrow S'$ as an alternative notation for  $\delta(S,(i,j))=S'$,
  and denote the initial state by $id$ instead of $Id$ for clarity.

\begin{prop}\label{dis} Let $M$ be a gasket automaton. For any $\bx,\by,\bz\in\Sigma^\infty$ we have
\begin{equation}\label{dis1}
\min\{T_M(\bx,\by), T_M(\bx,\bz)\}\le T_M(\by,\bz)+1.
\end{equation}
\end{prop}
\begin{proof} We will denote $T:=T_M$ for simplicity.
Clearly, \eqref{dis1} holds if either $T(\by,\bz)=\infty$ or any two of $\bx,\by,\bz$ are equal.
 So  we assume that $T(\by,\bz)<\infty$ and $\bx,\by,\bz$ are distinct.

Denote $\ell:=|\by\wedge\bz|$.  Then $T(\by,\bz)\ge\ell$, and at least one of $|\bx\wedge\by|\le\ell$ and  $|\bx\wedge\bz|\le\ell$ holds. Without loss of generality, assume that
$$
k=|\bx\wedge\by|\le |\bx\wedge\bz|.
$$
Suppose on the contrary that \eqref{dis1} is false, then $T(\bx,\by)>\ell+1$ and $T(\bx,\bz)>\ell+1$.

\textit{Case 1.}  $k<\ell$.

In this case we have
 $|\bx\wedge\by|=|\bx\wedge\bz|=k$,
 which together with $T(\bx,\by)>\ell+1$ imply that the first $(\ell+3)$-states (including $id$) of the itinerary of $(\bx,\by)$ are
\begin{equation}\label{dis2}
id\to(Id)^k\to (S_{uv})^{\ell+2-k}, \quad \text{ where }u,v\in\{\alpha,\beta,\gamma\}.
\end{equation}
So $(x_{k+1}, y_{k+1})\in P_{uv}$ and
$$
(\sigma^k(\bx),\sigma^k(\by))=(x_{k+1}v^{\ell+1-k}\cdots, y_{k+1}u^{\ell+1-k}\cdots).
$$
By $|\by\wedge\bz|=\ell$,
the first $\ell+1$ states (including $id$) of the itinerary of $(\bx,\bz)$ are the same as that of $(\bx,\by)$; in particular,  the $(\ell+1)$-th state is $S_{uv}$. Moreover,
 since $T(\bx,\bz)>\ell+1$, a prefix  of the itinerary of $(\bx,\bz)$  is also given by \eqref{dis2}.
 It follows that
 $\sigma^k(\bz)=y_{k+1}u^{\ell+1-k}\cdots$ and   $|\by\wedge\bz|\geq \ell+2$, which is a contradiction.
   Hence \eqref{dis1} holds in this case.

\textit{Case 2.} $k=\ell=|\bx\wedge\bz|$.

Then $x_1\dots x_k=y_1\dots y_k=z_1\dots z_k$ and $x_{k+1},y_{k+1}$ and $z_{k+1}$ are distinct. So the   itinerary of $(\bx,\by)$ is initialled by
\begin{equation}\label{dis3}
id \to(Id)^k\to (S_{uv})^2, \text{ where } u,v\in\{\alpha,\beta,\gamma\},
\end{equation}
 and    the itinerary of $(\bx,\bz)$ is initialled by
\begin{equation}\label{dis4}
id\to(Id)^k\to (S_{wv'})^2,\text{ where } w,v'\in\{\alpha,\beta,\gamma\}.
\end{equation}
The $(k+2)$-th transitions of \eqref{dis3} and \eqref{dis4} imply  $(x_{k+2},y_{k+2})=(v,u)$
and  $(x_{k+2}, z_{k+2})$ $=(v', w)$, and it follows that  $v=v'$. The $(k+1)$-th transitions imply
    $x_{k+1}\lhd_{uv}y_{k+1}$ and $x_{k+1}\lhd_{wv}z_{k+1}$, then $w\neq u$ by the uniqueness property, and  $y_{k+1}\lhd_{wu}z_{k+1}$ by the gathering condition.

Let $p$ be the largest integer such that
$$
(\sigma^k(\bx),\sigma^k(\by),\sigma^k(\bz))=(x_{k+1}v^p\cdots, y_{k+1}u^p\cdots, z_{k+1}w^p\cdots).
$$
Then all of $T(\bx,\by),$$ T(\by,\bz), $$T(\bx,\bz)$ are no less than $k+p+1$, and two of them are no larger than $k+p+2$ since $(x_{k+p+2},y_{k+p+2},z_{k+p+2}))\neq (v,u,w)$. The lemma holds in this case.

\textit{Case 3.} $k=\ell<|\bx\wedge\bz|$.

Since $x_{k+1}\neq y_{k+1}$, equation \eqref{dis3} still holds. Let $p$ be the largest integer such that
$$
(\sigma^k(\bx),\sigma^k(\by))=(x_{k+1}v^p\cdots, y_{k+1} u^p\cdots).
$$
Then $T(\bx,\by)=k+1+p$.

Let $q$ be the largest integer such that $\sigma^k(\bz)=x_{k+1}v^q\cdots$.
If $q\geq p-1$, then $T(\by,\bz)\geq k+p$ since $|\bx\wedge \bz|\geq k+p$, and the lemma holds in this case.
If $q\leq p-2$, then $T(\by,\bz)=k+q+1$ and
$$(\sigma^{k+q+1}(\bx), \sigma^{k+q+1}(\bz))=(v^2\cdots, \tilde v \eta \cdots),  \quad \text{ where } \tilde v\neq v.$$
Suppose $Id$ is not translated to $Exit$ by $(v,\tilde v)$, say, $Id\stackrel{(v,\tilde v)}\longrightarrow S$. If $S=S_{wv}$ for some $w\in\{\alpha,\beta,\gamma\}$, then $\tilde v\lhd_{vw} v$,
which violates the boundary condition.
 Thus $S\stackrel{(v,\eta)}\longrightarrow Exit$ and $T(\bx,\bz)\leq k+q+2$. The lemma  holds in this scenario.
This finishes the proof of \eqref{dis1}.
\end{proof}

 Now we can prove Theorem \ref{lem:feasible}.

 \begin{proof}[\textbf{Proof of Theorem \ref{lem:feasible}}]
Let $\bx,\by,\bz\in\Sigma^\infty$. Denote $\rho:=\rho_{M,\xi}$.

First, it is obvious that $T_M(\bx,\bx)=\infty$, so $\rho(\bx,\bx)=0$.

Secondly, $T_M(\bx,\by)=T_M(\by,\bx)$  by Lemma \ref{lem:symmetry}, so $\rho(\bx,\by)=\rho(\by,\bx)$.

Thirdly, by Proposition \ref{dis}, we have
$$
\rho(\bx,\by)+\rho(\by,\bz)=\xi^{T_M(\bx,\by)}+\xi^{T_M(\by,\bz)}
\geq \xi^{T_M(\bx,\bz)+1}=\xi\rho(\bx,\bz).
$$
The theorem is proved.
\end{proof}

Hence a gasket automaton $M$ induces a psuedo metric space, which we denote by $({\mathcal A_M}, \rho_M)$.
By \eqref{eq:xi2} and \eqref{eq:xi3}, we have
\begin{equation}\label{eq:xi22}
\rho_M(\bx,\by)\leq \xi^{-2}\rho_M([\bx],[\by]) \text{ and }
\end{equation}
\begin{equation}\label{eq:xi33}
\rho_M([\bx],[\bz])\leq \xi^{-3} (\rho_M([\bx],[\by])+\rho_M([\by],[\bz])).
\end{equation}

\section{\textbf{Topology automaton of fractal gasket}}\label{sec:structure}

 The neighbor graph (automaton) of self-similar sets is an important tool in fractal geometry, see \cite{BanM09,YZ18,RaoZhu16}.  The topology automaton
 is a simplified version of the neighbor automaton.

\begin{defn}[Topology automaton]\emph{Let $K$ be a fractal gasket generated by $\{\varphi_j\}_{j=1}^N$.
Let $\{\alpha,\beta,\gamma\}$ be a subset of $\Sigma\cup\{-1,-2,-3\}$ defined in Section 1.
Let $M_K$ be a triangle automaton satisfying:
For $i\neq j$,
$$
\delta(Id,(i,j))=\left \{
\begin{array}{ll}
S_{uv}, & \text{ if } u,v\in\Sigma \text{ and } \varphi_i(\omega_v)=\varphi_j(\omega_u),\\
Exit, &\text{ if } \varphi_i(K)\cap \varphi_j(K)=\emptyset.
\end{array}
\right .
$$
We call $M_K$ the \emph{topology automaton} of $K$.
}
\end{defn}

\begin{lem} The topology automaton of a fractal gasket is always a gasket automaton.
\end{lem}

\begin{proof} Let $K$ be a fractal gasket generated by $\{\varphi_j\}_{j=1}^N$.
 Since the functions $\varphi_i$ are all distinct,  $\varphi_i(K)$ can has at most one neighbor in $\theta$-direction for each $\theta\in\{\exp(2\pi \mathbf{i}k/6;~k=0,1,\dots, 5\}$.
 This verifies the uniqueness property of the gasket automaton.
 The  gathering condition and the boundary condition are obvious.
\end{proof}

For $A,B\subset\mathbb{R}^2$, let $\dist(A,B)=\min\{\|a-b\|;~a\in A, b\in B\}$.
Zhu and Yang \cite{YZ18}  defined  the sharp separation condition  for self-similar sets with uniform contraction ratio. We extend it to general self-similar sets.

\begin{defn}[Sharp separation condition]\emph{ A self-similar set $K$ is said to satisfy the \emph{sharp separation condition},
 if there exists a constant $C'>0$ such that for any $k\geq 1$ and   $I, J\in\Sigma^k$,  $\varphi_I(K)\cap\varphi_J(K)=\emptyset$ implies that
 $$dist(\varphi_I(K),\varphi_J(K))\ge C'\min\{\text{diam } \varphi_I(K), \text{diam } \varphi_J(K)\}.$$
 }
 \end{defn}

\begin{lem}\label{lem:sharp}  A fractal gasket   always satisfies the sharp separation condition.
\end{lem}
\begin{proof}
Let $K$ be a fractal gasket with IFS $\Phi=\{\varphi_i\}_{i=1}^N$.  Denote $K_i=\varphi_i(K)$ for $1\le i\le N$.
Let
$$C_1=\min\{dist(K_i,K_j); i,j\in \Sigma \text{ and } K_i\cap K_j=\emptyset\};$$
$$
C_2=\min\{dist(K_i,z); i\in \Sigma, z\in \{\omega_\alpha,\omega_\beta,\omega_\gamma\} \text{ and } z\not\in K_i\}.
$$

For any $I=x_1\dots x_k, J=y_1\dots y_k\in\Sigma^k$, suppose that $\varphi_I(K)\cap\varphi_J(K)=\emptyset$.
Let  $0\le\ell\le k-1$ be the largest integer such that
$\varphi_{x_1\dots x_\ell}(K)\cap\varphi_{y_1\dots y_\ell}(K)\ne\emptyset$.

If $x_1\dots x_{\ell}=y_1\dots y_{\ell}$, then
\begin{align*}
&dist(\varphi_I(K),\varphi_J(K))
= r_{x_1\dots x_{\ell}} dist(K_{x_{\ell+1}}, K_{y_{\ell+1}})\ge \frac{C_1}{\text{diam}(K)} \text{diam}(\varphi_I(K)).
\end{align*}

If $x_1\dots x_{\ell}\neq y_1\dots y_{\ell}$, let $z_0$ be the unique element  of $\varphi_{x_1\dots x_\ell}(K)\cap \varphi_{y_1\dots y_\ell}(K)$.
Let us assume that $z_0\not\in \varphi_{x_1\dots x_{\ell+1}}(K)$ without loss of generality.
Then
\begin{align*}
dist(\varphi_I(K),\varphi_J(K))&\ge dist(\varphi_{x_1\dots x_{\ell+1}}(K),z_0)\ge \frac{C_2}{\text{diam}(K)} \text{diam}(\varphi_I(K)).
\end{align*}
The lemma is proved.
\end{proof}

 Define $\pi:~\Sigma^\infty\to K$, which we call the \emph{coding map}, by
\begin{equation}
\big\{\pi(\bx)\big\}=\bigcap_{i\geq1} \varphi_{x_1\cdots x_i}(K).
\end{equation}
If $\pi(\bx)=x\in K$, then the sequence $\bx$ is called a \emph{coding} of $x$.

\begin{proof}[\textbf{Proof of Theorem \ref{thm:Holder}}]
Take $x,y\in K$. Let $\bx$ and $\by$ be a coding of $x$ and $y$, respectively.
Let $k=T(\bx,\by)$ be the surviving time in the topology automaton $M_K$. Then
the $k$-th cylinders containing $x$ and that containing $y$ either coincide or have non-empty intersection. It follows that
$$
\|x-y\|\leq 2(r^*)^{k}.
$$
On the other hand, the $(k+1)$-th cylinders containing $x$ and that containing $y$ are disjoint, so by the sharp separation condition, we have
$$
\|x-y\|\geq C'(r_*)^{k+1}
$$
where $C'$ is the constant in the sharp separation condition.
Recall that $s=\sqrt{\log r^*/\log r_*}$ and $\xi=(r_*)^s$.
So
$$\rho_{M_K, \xi}(\bx,\by)=\xi^k=(r_*)^{sk}=(r^*)^{k/s}.$$
Set $C=\max\{2, 1/(r_*C')\}$, we obtain the theorem.
\end{proof}

\section{\textbf{$\gamma$-isolated gasket automaton and simplification}}\label{sec:gamma}
In this section we   imposing additional conditions to the gasket automaton so that the automaton can
be simplified.

\begin{defn}[$\gamma$-isolated condition]\label{gamma-auto}
\emph{ Let $M$ be a gasket automaton.
We say $M$ satisfies the  \emph{$\gamma$-isolated condition} if  \\
\indent (i) $\{\alpha,\beta,\gamma\}\subset \Sigma$; \\
\indent (ii) The graph $(\Sigma,\mathcal{P}_{\alpha\gamma}\cup\mathcal{P}_{\beta\gamma})$ has no cycle;\\
 \indent (iii) $\gamma$ is isolated, in the sense that it is $\alpha\gamma$-isolated, $\beta\gamma$-isolated and $\alpha\beta$-isolated. }
\end{defn}

\begin{lem}If $K$ is a   fractal gasket such that
  $\{\omega_\alpha, \omega_\beta,\omega_\gamma\}\subset K$ and satisfying the top isolated condition,  then the topology automaton $M_K$  satisfies the  $\gamma$-isolated condition.
\end{lem}

\begin{proof} Let $K$ be a fractal gasket with IFS $\Phi=\{\varphi_i\}_{i=1}^N$.
If $i\lhd_{\alpha\gamma}j$, then $\varphi_j(\omega_\alpha)=\varphi_i(\omega_\gamma)$; denote by $c$ the contraction ratio of $\varphi_i$,  we have $$\varphi_j(\omega_\alpha)-\varphi_i(\omega_\alpha)=
\varphi_i(\omega_\gamma)-\varphi_i(\omega_\alpha)=c(\omega_\gamma-\omega_\alpha),$$
and it follows that  $\varphi_j(0,0)$ has larger second coordinate than that of $\varphi_i(0,0)$.
The same conclusion holds if $i\lhd_{\beta\gamma}j$.  This verifies (ii) in
Definition \ref{gamma-auto}. Clearly the top isolated condition implies that $\gamma$ is isolated.
\end{proof}


In the rest of this section, we always assume that $M$ is a  gasket automaton satisfying the $\gamma$-isolated
condition
such that
\begin{equation}
\mathcal{P}_{\alpha\gamma}\cup\mathcal{P}_{\beta\gamma}\ne\emptyset.
\end{equation}

For any $b\in\Sigma$, we say $b$ is \emph{double-maximal} in $M$ if $b$ is both $\alpha\gamma$-maximal and $\beta\gamma$-maximal. 


\begin{lem}\label{haschain}
There exists $(\tau,\kappa)\in\mathcal{P}_{\alpha\gamma}\cup\mathcal{P}_{\beta\gamma}$ such that $\kappa$ is double-maximal.
\end{lem}
\begin{proof}
Since $\mathcal{P}_{\alpha\gamma}\cup\mathcal{P}_{\beta\gamma}\ne\emptyset$, there exists $(b_1,b_2)\in\mathcal{P}_{\alpha\gamma}\cup\mathcal{P}_{\beta\gamma}$.
Let $(a_1,\dots, a_k)$ be the chain in $\mathcal{P}_{\alpha\gamma}\cup\mathcal{P}_{\beta\gamma}$
containing $(b_1,b_2)$, then $(\tau, \kappa)=(a_{k-1},a_k)$ is the desired edge.
\end{proof}

From now on,  we fix a pair $(\tau, \kappa)$ satisfying Lemma \ref{haschain}.
Moreover, we  assume that  $(\tau,\kappa)\in\mathcal{P}_{\alpha\gamma}$ without loss of generality.

If $\kappa$ has no $\alpha\beta$-predecessor, we set
\begin{equation}\label{fig_break1}
\mathcal{P}'_{\alpha\beta}=\mathcal{P}_{\alpha\beta}, \
\mathcal{P}'_{\alpha\gamma}=\mathcal{P}_{\alpha\gamma}\setminus\{(\tau,\kappa)\}, \text{ and }
\mathcal{P}'_{\beta\gamma}=\mathcal{P}_{\beta\gamma}.
\end{equation}
If $\kappa$ has a $\alpha\beta$-predecessor, we denote it by $\lambda$ and set
\begin{equation}\label{fig_break2}
\mathcal{P}'_{\alpha\beta}=\mathcal{P}_{\alpha\beta}, \
\mathcal{P}'_{\alpha\gamma}=\mathcal{P}_{\alpha\gamma}\setminus\{(\tau,\kappa)\}, \text{ and }
\mathcal{P}'_{\beta\gamma}=\mathcal{P}_{\beta\gamma}\setminus\{(\tau,\lambda)\}.
\end{equation}

Let $M'$ be the triangle automaton determined by $\mathcal{P}'_{\alpha\beta}, \mathcal{P}'_{\alpha\gamma}$ and $\mathcal{P}'_{\beta\gamma}$, and we call it a \emph{one-step simplification} of $M$.
If \eqref{fig_break1} holds, we call $M'$ a $(\tau, \kappa)$-simplification, otherwise, we call
$M'$ a $(\tau,\kappa,\lambda)$-simplification.

\begin{figure}[H]
\centering
\includegraphics[width=14 cm]{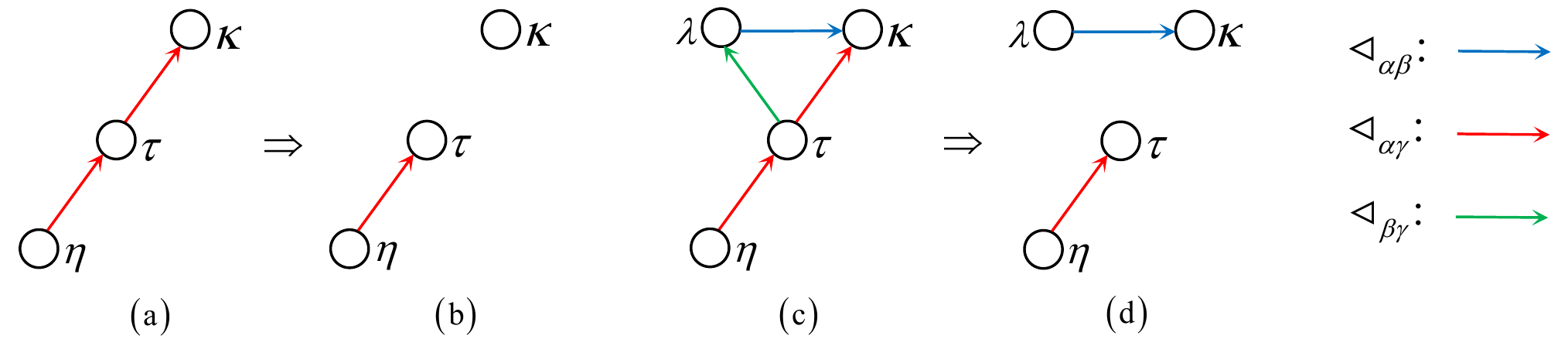}\\
\caption{(a)$\Rightarrow$(b) illustrates a $(\tau,\kappa)$-simplification, while  (c)$\Rightarrow$(d) illustrates a $(\tau,\kappa,\lambda)$-simplification.}\label{break}
\end{figure}

\begin{lem}\label{M'gasketauto}
Let $M'$ be the one-step simplification of $M$, then $M'$ is also a  gasket automaton satisfying the $\gamma$-isolated condition.
\end{lem}
\begin{proof}
Notice that $(\Sigma,\mathcal{P}'_{\alpha\beta})=(\Sigma,\mathcal{P}_{\alpha\beta})$, moreover, $(\Sigma,\mathcal{P}'_{\alpha\gamma})$ and $(\Sigma,\mathcal{P}'_{\beta\gamma})$ are subgraphs of $(\Sigma,\mathcal{P}_{\alpha\gamma})$ and $(\Sigma,\mathcal{P}_{\beta\gamma})$ respectively, so $M'$ satisfies the unique property  and the boundary condition in the definition of gasket automaton.
Also, item (ii) in Definition \ref{gamma-auto} holds.

If three edges with vertices $a,b,c$ satisfies the assumptions \ding{172}\ding{173}\ding{174} in gathering condition, then we call the set
$\{(a,b),(b,c),(c,a)\}$ a \emph{family} of $M$. Since one edge in a family can determine the other two edges,
we obtain that any two families are edge-disjoint.

If $M'$ is a $(\tau,\kappa)$-simplification, then $(\tau,\kappa)$ does not belong to any family, so the simplification does not affect any family, and $M'$ satisfies the gathering condition. If $M'$ is a $(\tau,\kappa,\lambda)$-simplification, then the simplification deletes two members of a family, so $M'$ still
satisfies the gathering condition. The lemma is proved.
\end{proof}

\begin{figure}[H]
\centering
\includegraphics[width=0.35 \textwidth]{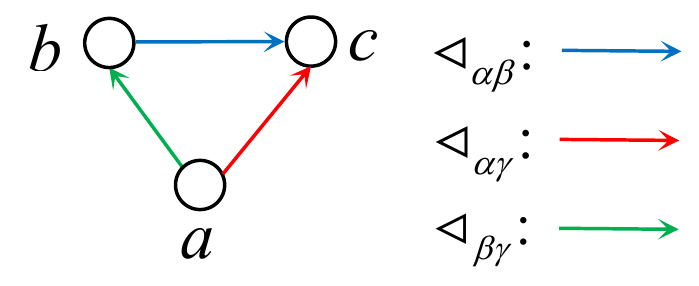}
\caption{A family in $M$.}\label{three}
\end{figure}

\begin{lem}\label{max_element}
Let $M'$ be the one-step simplification of $M$. Then\\
\indent \emph{(i)} both $\tau$ and $\kappa$ are double-maximal in $M'$, and $\kappa$ is $\alpha\gamma$-isolated in $M'$.\\
\indent \emph{(ii)} $\kappa\notin\{\alpha,\gamma\}$.
\end{lem}
\begin{proof}
(i)   That $\kappa$ is double-maximal in $M'$  since it is double-maximal in $M$. The one-step simplification
 deletes the edges connecting $\tau$ to its $\alpha\gamma$-successor and $\beta\gamma$-successor (if exists),
 so the other assertions in (i) hold. (See Figure \ref{break} (b) or (d) for an illustration.)

(ii) Since $\tau\lhd_{\alpha\gamma} \kappa$,  that $\alpha$ is $\alpha\gamma$-minimal
and $\gamma$ is $\alpha\gamma$-isolated imply  $\kappa\notin\{\alpha, \gamma\}$.
\end{proof}

The following theorem plays a crucial role in this paper.

\begin{thm}\label{main}
Let $M'$ be the one-step simplification of $M$. Then for any $\bx,\by\in\Omega:=\{\omega\kappa^{\infty};\omega\in\Sigma^*\}$ there exists a bijection $g:\Omega\to \Omega$ such that
\begin{equation}\label{eq-key}
|T_M(\bx,\by)-T_{M'}(g(\bx),g(\by))|\le 5.
\end{equation}
\end{thm}

In Section \ref{sec:mapg}, we give the construction of the map $g$,  and  we   prove Theorem \ref{main}  in Section \ref{sec:time}.

\section{\textbf{Proof of Theorem \ref{spaceLip} and Theorem \ref{thm:Lip}}}\label{sec:proof}

 Let $M$ be a  gasket automaton satisfying the $\gamma$-isolated condition.

We denote by $M^*$  the gasket automaton  determined by $\mathcal{P}^*_{\alpha\gamma}=\mathcal{P}^*_{\beta\gamma}=\emptyset$ and $\mathcal{P}^*_{\alpha\beta}=\mathcal{P}_{\alpha\beta}$, and  we call it the \emph{final-simplification} of $M$.

\begin{proof}[\textbf{Proof of Theorem \ref{spaceLip}}]   Recall that $\Omega=\{\omega\kappa^{\infty};\omega\in\Sigma^*\}.$
Let $g:\Omega\rightarrow\Omega$ be the map given in Theorem \ref{main}. Since for any $\bx,\by\in\Omega$,
$$
T_M(\bx,\by)-5\le T_{M'}(g(\bx),g(\by))\le T_M(\bx,\by)+5,
$$
which implies that
$\xi^5\rho_M(\bx,\by)\le\rho_{M'}(g(\bx),g(\by))\le\xi^{-5}\rho_M(\bx,\by).$
Hence $g$ is  bi-Lipschitz.

Define $[\bx]_M=\{\by\in \Sigma^\infty;~ \rho_M(\bx,\by)=0\}$ and set $\Omega_M=\{[\bx]_M; ~\bx\in \Omega\}$; similarly we define $[\bx]_{M'}$ and  $\Omega_{M'}$.
Define $\widetilde g: \Omega_M\to \Omega_{M'}$ by $\widetilde g([\bx]_M)=[g(\bx)]_{M'}$. We claim that
$\widetilde g$ is bi-Lipschitz. First,
$$
\begin{array}{rl}
& \rho_{M'}(\widetilde g([\bx]_M), \widetilde g([\by]_M))=\rho_{M'}([g(\bx)]_{M'}, [g(\by)]_{M'})\\
\leq &\rho_{M'}(g(\bx), g(\by))\leq \xi^{-5} \rho_M(\bx,\by)\\
\leq & (\xi^{-5})(\xi^{-2}) \rho_M([\bx]_M, [\by]_M).  \quad \text{(By \eqref{eq:xi22}.)}
\end{array}
$$
In the same manner, we can prove the other direction inequality. The claim is proved.

Moreover, by Lemma \ref{Omegadense}, $\Omega_M$
is dense in  $({\mathcal A}_M,\rho_M)$ and $\Omega_{M'}$ is dense in $({\mathcal A}_{M'},\rho_{M'})$. So
${\mathcal A}_M\simeq {\mathcal A}_{M'}$  by Lemma \ref{extendLip}.
\end{proof}

Using  Theorem \ref{spaceLip} repeatedly, we obtain

\begin{coro}\label{finalspaceLip}
Let $M$ be a gasket automaton and $M^*$ be the final-simplification of $M$. Then $(\mathcal{A}_M,\rho_M)\simeq(\mathcal{A}_{M^*},\rho_{M^*})$.
\end{coro}

\begin{proof} Notice that a gasket automaton $M$ admits a one-step simplification provided $\mathcal{P}_{\alpha\gamma}\cup\mathcal{P}_{\beta\gamma}\ne\emptyset$ (see Section 5). Then there exists a sequence
$$
M=M_0,\ M_1,\ \dots,\ M_q=M^*
$$
such that $M_{j+1}$ is the one-step simplification of $M_j$ for each $0\le j\le q-1$.
So the result is a consequence of Theorem \ref{spaceLip}.
\end{proof}

 Let $E$ and $F$ be the fractal gaskets defined in Theorem \ref{thm:Lip}, and let
 $M_E$ and $M_F$ be the topology automata of $E$ and $F$, respectively.
  Without loss of generality, we assume that $\alpha=1$ and $\beta=2$ for both $E$ and $F$.

If $\omega_\gamma\in E$, we set  $M^*_E=(M_E)^*$ to be  the final-simplification of $M_E$, otherwise, we set  $M_E^*=M_E.$
Similarly we define $M^*_F$.

\begin{remark}\emph{  Actually, if $\omega_\gamma\not\in E$, then the topology automaton $M_E$ only records the horizontal connective relations among $\varphi_j(E), j\in\Sigma$; more precisely, for any $i,j\in\Sigma$, $\varphi_i(E)\cap\varphi_j(E)\ne\emptyset$ if and only if $i\lhd_{\alpha\beta}j$ or $j\lhd_{\alpha\beta}i$. Thus $\mathcal{P}_{\alpha\gamma}\cup\mathcal{P}_{\beta\gamma}=\emptyset$.
}
\end{remark}

\begin{lem}\label{finalEFLip}
There exists an isometry  $f:(\mathcal{A}_{M^*_E},\rho_{M^*_E})\to(\mathcal{A}_{M^*_F},\rho_{M^*_F})$.
\end{lem}
\begin{proof}
Let $I=\{a_1,a_2,\dots,a_k\}\subset\Sigma$ be a horizontal-block of $E$. By the definition of horizontal-block and neighbor automaton $M_E$, we have
$$
a_j\lhd_{\alpha\beta}a_{j+1}\text{ in }M_E,\ 1\le j\le k-1.
$$
By  assumptions of Theorem \ref{thm:Lip}, there is a size-preserving bijection from the collection of horizontal-blocks of $E$ to that of $F$, which we denote by $\widehat{h}$. That is,
\begin{equation}\label{finalEFLip1}
\widehat{h}(I)=\{b_1,b_2,\dots,b_k\}
\end{equation}
is a horizontal-block of $F$.
Define $h:\Sigma\to\Sigma$ by $h(a_j)=b_j$, that is, if $a_j$ is the $j$-th element of a horizontal block $I$ of $E$, then we define $h(a_j)$ to be the $j$-th element of $\widehat{h}(I)$. Then, for any $r,s\in\Sigma$,
\begin{equation}\label{finalEFLip2}
r\lhd_{\alpha\beta}s\text{ in }M_E^*\quad\text{if and only if}\quad h(r)\lhd_{\alpha\beta}h(s)\text{ in }M_F^*.
\end{equation}

Now we define $f:\Sigma^\infty\to\Sigma^\infty$ by
$
f((x_i)_{i=1}^\infty)=(h(x_i))_{i=1}^\infty.
$
Clearly, $f$ is a bijection and   for any $\bx,\by\in\Sigma^\infty$,
 $T_{M_E^*}(\bx,\by)=T_{M_F^*}(f(\bx),f(\by)).$
It follows that  $[\bx]\mapsto [f(\bx)]$ is an isometry from  $\mathcal{A}_{M^*_E}$ to $\mathcal{A}_{M^*_F}$.
\end{proof}

\begin{proof}[\textbf{Proof of Theorem \ref{thm:Lip}}] Fix $\xi\in (0,1)$. We have
$$
\mathcal{A}_{M_E}\simeq \mathcal{A}_{M^*_E} \simeq \mathcal{A}_{M^*_F} \simeq \mathcal{A}_{M_F},
$$
where the first and third relations are due to Corollary \ref{finalspaceLip},
and the secondly relation is by Lemma \ref{finalEFLip}.
Next, by Theorem \ref{thm:Holder},  we have that
$E$ is bi-H\"{o}lder equivalent to $\mathcal{A}_{M_E}$ and
$F$ is bi-H\"{o}lder equivalent to $\mathcal{A}_{M_F}.$
Thus $E$ is bi-H\"{o}lder equivalent to $F$. (Here we use the fact that $({\mathcal A_M}, \rho_{M,\xi})$
is bi-H\"older equivalent to $({\mathcal A_M}, \rho_{M,\xi'})$ provided $\xi'\in (0,1)$.)

Especially, if both $E$ and $F$ have uniform contraction ratio $r$, then setting $\xi=r$, we obtain $E\simeq F$.
\end{proof}

\section{\textbf{A universal map from $\Omega$ to $\Omega$}}\label{sec:mapg}

Let $\Sigma=\{1,2,\dots,N\}$ with $N\ge 4$. Let $\tau,\kappa,\alpha,\gamma\in\Sigma$   be distinct except that
 $\tau=\alpha$ is allowed.  Set \begin{equation}\label{Omega}
\Omega=\{\omega\kappa^{\infty};\omega\in\Sigma^*\}.
 \end{equation}
  In this section, we construct a bijection $g:\Omega\to\Omega$, which can be realized by a transducer (see Appendix B).  In next section, we will show that the map $g$ is the desired map in Theorem \ref{main}.

 \begin{remark} \emph{(i) The discussion of this section is purely symbolic: it is irrelevant to metric or automaton. }

 \emph{(ii) If $\tau$ and $\kappa$ come from a one-step simplification, then
  $\kappa\not\in\{ \alpha, \gamma\}$ by Lemma \ref{max_element},  and $\tau\neq \gamma$ since $\gamma$ is isolated.}
 \end{remark}

\subsection{Segment decomposition}
 \ \\
\indent First we introduce two decompositions of sequences in $\Omega$.
Set
\begin{equation}\label{CM}
\mathcal{C}_M:=\{\tau\gamma^k;k\ge 2\}\cup\{\kappa\alpha^{k}\kappa\gamma;k\ge 0\}.
\end{equation}

\begin{defn}[$M$-decomposition]\label{M-segment}
\emph{Let $\bx=(x_i)_{i=1}^\infty\in\Omega$.  The longest prefix  $X_1$ of $\bx$ satisfying
$X_1\in {\mathcal C}_M\cup \Sigma$ is called the \emph{$M$-initial segment} of $\bx$.}

\emph{Inductively, each $\bx=(x_i)_{i=1}^{\infty}\in\Omega$ can be uniquely written as
$\bx=\prod_{j=1}^\infty X_j:=X_1X_2\cdots X_k \cdots,$
where $X_k$ is the $M$-initial segment of $\prod_{j\geq k} X_j$. We call $(X_j)_{j\geq 1}$ the \emph{$M$-decomposition} of $\bx$.}
\end{defn}

Next we define  $M'$-decomposition. Set
\begin{equation}\label{CM'}
\mathcal{C}_{M'}=\{\kappa\alpha^k\kappa\gamma;k\ge 0\}\cup\{\kappa\alpha^{k}\kappa\gamma\gamma;k\ge 0\}\cup\{\tau\gamma\gamma\},
\end{equation}

\begin{defn}[$M'$-decomposition]\label{M'-segment}
\emph{Let $\bu=(u_i)_{i=1}^\infty\in\Omega$.  A word $U_1$  is called the $M'$-initial segment of $\bu$, if
it is the longest prefix of $\bu$ such that $U_1\in {\mathcal C}_{M'}\cup \Sigma$. Similar as above,   we define the
\emph{$M'$-decomposition} of $\bu$.
 }
\end{defn}

Two words are said to be \emph{comparable}, if one  is a prefix of the other.

\begin{remark}\label{rem:observe} \emph{Here are two useful observations.}

\emph{ (i) If two elements in ${\mathcal C}_M$ are  comparable,
then both of them are of the form $\tau\gamma^k$. If two elements in ${\mathcal C}_{M'}$ are  comparable,
then one of them is $\kappa\alpha^k\kappa\gamma$ and another one is $\kappa\alpha^k\kappa\gamma\gamma$.
}

\emph{(ii) Let $W\in {\mathcal C}_M\cup {\mathcal C}_{M'}$. Then $W$ is initialled by a word in
$\{\kappa \alpha, \kappa\kappa, \tau\gamma\}$. Moreover, these words cannot  appear in  $W$ except as a prefix.
}
\end{remark}

\subsection{Construction of $g$}

First we define $g_0:\mathcal{C}_M\cup\Sigma\rightarrow\mathcal{C}_{M'}\cup\Sigma$ by
$$
g_0:\left\{\begin{array}{rl}
\tau\gamma^k&\mapsto\kappa\alpha^{k-2}\kappa\gamma,\ k\ge 2;\\
\kappa\alpha^k\kappa\gamma&\mapsto\kappa\alpha^{k-1}\kappa\gamma\gamma,\ k\ge 1;\\
\kappa\kappa\gamma&\mapsto\tau\gamma\gamma;\\
i&\mapsto i, \forall i\in\Sigma.
\end{array}
\right.
$$

Clearly $g_0:\mathcal{C}_M\cup\Sigma\rightarrow\mathcal{C}_{M'}\cup\Sigma$ is a bijection.  We define $g:\Omega\rightarrow\Omega$ by
\begin{equation}\label{gxdecomposition}
g(\bx)=\prod_{j=1}^\infty g_0(X_j),
\end{equation}
where $(X_j)_{j=1}^\infty$ is the $M$-decomposition of $\bx$.
(Notice that  any $\bx\in \Omega$ has
a $M$-decomposition of the form $(X_j)_{j=1}^\ell(\kappa)^\infty$, so $g(\bx)=(\prod_{j=1}^\ell g_0(X_j))(\kappa)^\infty\in\Omega$.)

\begin{remark}\emph{Rao and Zhu \cite{RaoZhu16} constructed the first map of this type related to two neighbor automata of fractal squares,
 based on geometrical observations.  The  map $g$ above is an improvement of the map in \cite{RaoZhu16}.
 }
\end{remark}

\begin{prop}\label{coincide}
Let $\bx=x_1x_2\dots, \bu=u_1u_2\dots=g(\bx)$.\\
\indent \emph{(i)} If $(X_j)_{j\geq 1}$ is the $M$-decomposition of $\bx$, then the $M'$-decomposition of $g(\bx)$ is $\left (g_0(X_j)\right )_{j\geq 1}$.\\
\indent \emph{(ii)} Similarly, if $(U_j)_{j\geq 1}$ is the $M'$-decomposition of $\bu$, then the $M$-decomposition of $h(\bu)$ is $\left ( g^{-1}_0(U_j) \right )_{j\geq 1}$, where $h(\bu)=\prod_{j=1}^\infty g_0^{-1}(U_j)$.\\
\indent \emph{(iii)} The map $g:\Omega\rightarrow\Omega$ is a bijection.
\end{prop}

\begin{proof} We denote $U\lhd W$ if  $U$ is a prefix of $W$.

(i) To prove the first statement, we only need to show that $U_1=g_0(X_1)$.

Suppose on the contrary that $U_1\neq g_0(X_1)$.
First, $U_1\lhd g_0(X_1)$ is impossible, since $g_0(X_1)\in {\mathcal C}_{M'}$ and we always choose
the longest one to be the initial segment.
Let $p$ be the least integer such that  $U_1\lhd g_0(X_1)\cdots g_0(X_p)$, then $p\geq 2$.

If $|X_j|=1$ for all $j=1,\dots, p$, then $\bx$ is initialled by $x_1\dots x_{|U_1|}=U_1$, which forces that $U_1\in {\mathcal C}_{M'}\setminus {\mathcal C}_M$, so $U_1=\kappa \alpha^k\kappa\gamma\gamma$ is the only choice. But then the initial segment of $\bx$ should be $\kappa\alpha^{k}\kappa\gamma$, a contradiction.

So at least one $X_j\in {\mathcal C}_M$. Suppose it happens for $j_0$. Since elements in ${\mathcal C}_M$
are initialled by $\kappa\kappa$, $\kappa\alpha$ or $\tau\gamma$ , we deduce that $j_0=1$.
Thus $g_0(X_1)$ is a proper prefix of $U_1$, which forces that $g_0(X_1)=\kappa\alpha^k\kappa\gamma$
and $U_1=\kappa \alpha^k\kappa\gamma\gamma$. But then $X_1=\tau \gamma^{k+2}$ and $x_{k+4}=\gamma$, which contradicts that $X_1$ is an initial segment.
Item (i) is proved.

(ii) For the second assertion, we only need to show that $X_1=g_0^{-1}(U_1)$, which is almost the same as the proof of (i).

Suppose on the contrary that $X_1\neq g_0^{-1}(U_1)$.
Then $X_1\lhd g_0^{-1}(U_1)$ is impossible since $g_0^{-1}(U_1)\in {\mathcal C}_{M}$.
Let $p\geq 2$ be the least integer such that $X_1\lhd g_0^{-1}(U_1)\cdots g_0^{-1}(U_p)$.

If $|U_j|=1$ for all $j=1,\dots, p$, then $\bu$ is initialled by $X_1$,  so $X_1\in
{\mathcal C}_M\setminus {\mathcal C}_{M'}$
and $X_1=\tau\gamma^k (k\ge 3)$ is the only choice. But  then the initial segment of $\bu$ should be $\tau\gamma\gamma$, a contradiction.

By the same reason as item (i),  we have $U_{1}\in {\mathcal C}_{M'}$. Thus $g_0^{-1}(U_1)$ is a proper prefix of $X_1$, which forces that  $g_0^{-1}(U_1)=\tau\gamma^k (k\geq 2)$
and $X_1=\tau\gamma^\ell (\ell>k)$. But then $U_1=\kappa\alpha^{k-2}\kappa\gamma$ and $u_{k+2}=\gamma$,
 which is  a contradiction.
Item (ii) is proved.

(iii)   From (i) and (ii) we have $h\circ g=g\circ h=id$, so $g$ is a bijection.
\end{proof}

Let $\sigma:\Sigma^\infty\to\Sigma^\infty$ be the shift operator defined by $\sigma((x_k)_{k\ge 1})=(x_k)_{k\geq 2}.$

\begin{lem}\label{lem-length}
Let $\bx=(x_k)_{k\geq 1}, \by=(y_k)_{k\geq 1}\in\Omega$. Then
$$|g(\bx)\wedge g(\by)|\geq |\bx\wedge\by|-2.$$
In other words,  $u_1\cdots u_k$ is determined by $x_1\cdots x_{k+2}$, where $k\ge 1$.
\end{lem}
\begin{proof}  Denote $\bu=g(\bx), \bv=g(\by)$.
Let $(X_j)_{j=1}^\infty$ and $(Y_j)_{j=1}^\infty$ be the $M$-decompositions of $\bx$ and $\by$ respectively.
Denote $U_1=g_0(X_1)$ and $V_1=g_0(Y_1)$.
Let $k=|\bx\wedge \by|-2\geq 1$.

We first prove the lemma in case of $X_1\neq Y_1$. Without loss
of generality, we assume that $|Y_1|\geq |X_1|$.

\ding{172} Suppose $|X_1|=1$. In this case we must have $Y_1=\kappa\alpha^q\kappa\gamma$ and $q\geq k$.
Moreover,  we have $|X_i|=1$ as long as $i\leq |Y_1|-2$, since a word in ${\mathcal C}_M$ is initialled  by
$\kappa\kappa$, $\kappa\alpha$ or $\tau\gamma$.
Since, $k+1\leq |Y_1|-2$, this implies that  $x_1\dots x_{k+1}$ (which is the same as $y_1\dots y_{k+1})$ is a prefix of $\bu$.
Since $V_1=\kappa\alpha^{q-1}\kappa\gamma\gamma$ is initialled by $y_1\dots y_q$, we have
$|\bu\wedge\bv|\geq |y_1\dots y_{k+1}\wedge y_1\dots y_{q}|=\min\{k+1,q\}\geq k.$

\ding{173} Suppose $|X_1|>1$.
If $X_1=\tau\gamma^{\ell}$,  then    $Y_1=\tau \gamma^{q}$ with $q>\ell$, so  $|\bu\wedge\bv|=|U_1\wedge V_1|=|\kappa\alpha^{\ell-2}\kappa\gamma\wedge \kappa\alpha^{q-2}\kappa\gamma|= \ell-1=k$.
If  $X_1=\kappa\alpha^{\ell}\kappa\gamma$, then  $Y_1=\kappa\alpha^q\kappa\gamma$ with $q>\ell$, so
$|\bu\wedge \bv|=|\kappa \alpha^{\ell-1}\kappa\gamma\wedge \kappa \alpha^{q-1}\kappa\gamma |= \ell=|\bx \wedge \by|-1 =k+1$.

Hence the lemma is valid if $X_1\neq Y_1$.

If $X_1=Y_1$, denote $p=|X_1|$. Denote $\ba=\sigma^p(\bx)$ and $\bb=\sigma^p(\by)$.
Hence the lemma holds for $\bx$ and $\by$ if and only if it holds for $\ba$ and $\bb$.
 So the lemma can be proved by induction.
\end{proof}

In Appendix B, we give a transducer which can realize the map $g$. The transducer provides an alternative proof of Lemma \ref{lem-length}.

\section{\textbf{Proof of Theorem \ref{main}}}\label{sec:time}

Let $M$ be a  gasket automaton satisfying the $\gamma$-isolated condition, and let $M'$ be a one-step simplification of $M$. For
$\bx, \by\in\Omega=\{\omega\kappa^{\infty};\omega\in\Sigma^*\}$, denote $\bu=g(\bx)~{\rm and}~ \bv=g(\by)$. Let $(X_j)_{j=1}^{\infty}$, $(Y_j)_{j=1}^{\infty}$ be the
$M$-decompositions of $\bx, \by$ respectively, and $(U_j)_{j=1}^{\infty}$, $(V_j)_{j=1}^{\infty}$ be the $M'$-decompositions of $\bu, \bv$ respectively. Clearly, we always have
\begin{equation}\label{eq:Less}
T_{M'}(\bx,\by)\leq T_M(\bx,\by).
\end{equation}

Recall that $\mathcal{C}_M=\{\tau\gamma^k;k\ge 2\}\cup\{\kappa\alpha^{k}\kappa\gamma;k\ge 0\}$ and
 $\mathcal{C}_{M'}=\{\kappa\alpha^k\kappa\gamma;k\ge 0\}\cup\{\kappa\alpha^{k}\kappa\gamma\gamma;k\ge 0\}\cup\{\tau\gamma\gamma\}$.
One should keep in mind that
\begin{equation}\label{eq:cycle}
S_{uv}\stackrel{(i,j)}\longrightarrow S_{uv} \text{ if and only if } (i,j)=(v,u)
\end{equation}
in both $M$ and $M'$. Since $\gamma$ is isolated,   if $\tilde \gamma\in \Sigma\setminus\{\gamma\}$,  then
\begin{equation}\label{eq:gammaExit}
Id\stackrel{(\gamma,\tilde \gamma)}\longrightarrow Exit \text{ and }
Id\stackrel{(\tilde \gamma,\gamma)}\longrightarrow Exit
\end{equation}
in both $M$ and $M'$.

\begin{lem}\label{gotoExit}
 Let $\ba,\bb\in\Omega$ such that $a_1\neq b_1$.  If $\ba=\alpha^k\kappa\gamma\cdots (k\geq 0)$,
 then
 $T_M(\ba,\bb)\leq 2.$
\end{lem}

\begin{proof} 
Suppose  $(a_1,b_1)\in\mathcal{P}_M$, then
$Id\stackrel{(a_1,b_1)}\longrightarrow S_{uv}$ for some $u,v\in\{\alpha,\beta,\gamma\}$.
To prove the lemma,  by \eqref{eq:cycle}, we only need to show that $a_2\neq v$ or $a_3\neq v$.

If $k=0$, then $(a_1,b_1)=(\kappa,b_1)$ and $a_2=\gamma$. Since $\kappa$ is double-maximal in both $M$, that is, $\kappa$ is $\alpha\gamma$-maximal and $\beta\gamma$-maximal,  we have $v\ne\gamma$. So $a_2\neq v$.

If $k=1$, then $a_2=\kappa$ and $a_3=\gamma$.
 Since $\kappa\neq \gamma$,  clearly either  $a_2\ne v$ or $a_3\neq v$.

If $k>1$, then $(a_1,b_1)=(\alpha,b_1)$ and  $a_2=\alpha$.
Since $\alpha$ is $\alpha\gamma$-minimal and $\alpha\beta$-minimal $M$,  we have $v\neq \alpha$.
The lemma is proved.
\end{proof}

\begin{lem}\label{case-id}
\emph{(i)} Let $\bx,\by\in\Omega$. If $x_1=y_1$ and $X_1\neq Y_1$, then
$$
T_M(\bx,\by)\leq |\bx\wedge\by|+2.
$$
\emph{(ii)} Let $\bu,\bv\in\Omega$. If $u_1= v_1$ and $U_1\neq V_1$, then
$$
T_{M'}(\bu,\bv)\leq |\bu\wedge\bv| +2.
$$
\end{lem}

\begin{proof}
(i) Let $k=|\bx\wedge\by|$, then $k\ge 1$ and
$x_{k+1}\ne y_{k+1}.$
By $X_1\neq Y_1$ we know that at least one of $X_1$ and $Y_1$ is in $\mathcal{C}_M$, say $X_1\in\mathcal{C}_M$. We consider two cases according to $X_1$.

\ding{172} $X_1=\tau\gamma^\ell(\ell\ge 2)$.  In this case, we have $k\le\ell+1$, for otherwise $Y_1=X_1$.

If $k\le\ell$, then $x_{k+1}=\gamma$, so $(x_{k+1},y_{k+1})=(\gamma,\widetilde{\gamma})$, where $\widetilde{\gamma}\in\Sigma\setminus\{\gamma\}$; if $k=\ell+1$, then $Y_1=\tau\gamma^s$ with $s>\ell$, which implies $(x_{k+1},y_{k+1})=(\widetilde{\gamma},\gamma)$.
Hence, by formula \eqref{eq:gammaExit},  $Id$ is transferred to $Exit$ by $(x_{k+1},y_{k+1})$,
so $T_M(\bx,\by)=k$.

\ding{173} $X_1=\kappa\alpha^\ell\kappa\gamma(\ell\ge 0)$. In this case, we have $k\le\ell+2$, for otherwise $X_1=Y_1$.

If $k\le\ell+1$, then $x_{k+1}x_{k+2}\cdots=\alpha^p\kappa\gamma\cdots (p\geq 0)$,
  and we have $T_M(\bx,\by)\leq k+2$ by Lemma \ref{gotoExit}.
If $k=\ell+2$, then $x_{k+1}=\gamma$ and $(x_{k+1},y_{k+1})=(\gamma,\widetilde{\gamma})$, so $T_M(\bx,\by)=k$.
This complete the proof of (i).

(ii) Using item (i) we have just proved, we have
$$
|\bu\wedge \bv|\leq T_{M'}(\bu,\bv)\leq T_M(\bu,\bv)\leq |\bu\wedge \bv|+2.
$$
The lemma is proved.
\end{proof}

\begin{lem}\label{tool1}
Let $\bx=x_1s^{k}\cdots$, where $k\ge 2$ and $s\in\{\alpha,\beta,\gamma\}$. Then
$$
g(\bx)=\left \{
\begin{array}{ll}
\kappa\alpha^{k-2}\cdots, &\text{ if } x_1s^{k}=\tau\gamma^{k} ;\\
x_1s^{k-2}\cdots, &\text{ otherwise.}
\end{array}
\right .
$$
\end{lem}
\begin{proof} By Lemma \ref{lem-length}, $u_1\cdots u_{k-1}$ is determined by $x_1\dots x_{k+1}$. So the lemma holds since
$g(\tau\gamma^{k}\kappa^\infty)=\kappa\alpha^{k-2}\kappa\gamma\kappa^\infty$ and
$g(x_1s^{k}\kappa^\infty)=x_1s^{k}\kappa^\infty$.
\end{proof}

\begin{lem}\label{lem:key-1} Let  $\bx,\by\in\Omega$. If  $X_1\ne Y_1$, then
\begin{equation}\label{eq-key-1}
T_M(\bx,\by)-T_{M'}(\bu,\bv)\le 5.
\end{equation}
\end{lem}

\begin{proof}
Let $S_{M,1}$ be the first state of the itinerary of $(\bx,\by)$ in $M$ (after the initial state $id$).
If  $S_{M,1}=Exit$,  obvious  \eqref{eq-key-1} holds.

If $S_{M,1}=Id$, then  $x_1=y_1$, so
$T_M(\bx,\by) \leq |\bx\wedge\by|+2$ by Lemma \ref{case-id}.
By Lemma \ref{lem-length}, we have  $T_{M'}(\bu,\bv)\geq |\bu\wedge \bv|\geq |\bx\wedge\by|-2$.
Hence  \eqref{eq-key-1} holds in this case.

Finally, we deal with the case $S_{M,1}\in Q\setminus\{Id, Exit\}$.
  Suppose $\tau$ has a $\beta\gamma$-successor $\lambda$ in $M$,  then $M'$ is the $(\tau,\kappa,\lambda)$-simplification of $M$ with
$\mathcal{P}'_{\alpha\beta}=\mathcal{P}_{\alpha\beta}, \ \mathcal{P}'_{\alpha\gamma}=\mathcal{P}_{\alpha\gamma}\setminus\{(\tau,\kappa)\} \text{ and } \mathcal{P}'_{\beta\gamma}=\mathcal{P}_{\beta\gamma}\setminus\{(\tau,\lambda)\}.$
(See Figure \ref{xy}.)

\begin{figure}[H]
\centering
\includegraphics[width=9 cm]{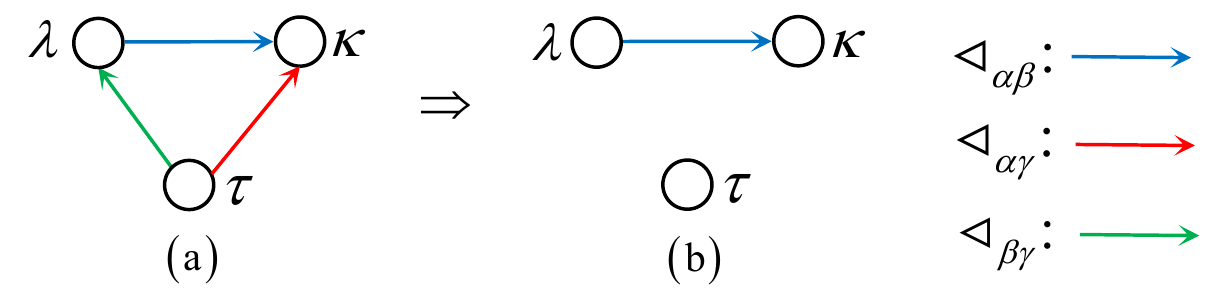}\\
\caption{$(\tau,\kappa,\lambda)$-simplification}\label{xy}
\end{figure}

 Denote $k=T_M(\bx,\by)$. If $k\le 5$,  \eqref{eq-key-1} holds trivially. Now we assume $k\ge 6$. Denote $S_{M,1}=S_{rs}$, where $r,s\in\{\alpha,\beta,\gamma\}$. Then the itinerary of $(\bx,\by)$ is $id\to (S_{rs})^{k}\to Exit$, so $(x_1,y_1)\in\mathcal{P}_{rs}$ in $M$ and
$$
\bx=x_1s^{k-1}\cdots,\quad \by=y_1r^{k-1}\cdots.
$$

\textit{Case 1. } $x_1s^{k-1},y_1r^{k-1}\ne\tau\gamma^{k-1}$.

By Lemma \ref{tool1} we have
$\bu=x_1s^{k-3}\cdots, \quad \bv=y_1r^{k-3}\cdots.$
 We claim that
\begin{equation}\label{lem-key-2eq2}
(x_1,y_1)\notin\{(\tau,\kappa),(\tau,\lambda)\}\cup\{ (\kappa,\tau),(\lambda,\tau)\}.
\end{equation}
If $(x_1,y_1)=(\tau,\kappa)$, then $x_1\lhd_{\alpha\gamma}y_1$, which implies that $s=\gamma$, which contradicts  $x_1s^{k-1}\neq \tau\gamma^{k-1}$. Similarly, we can eliminate the other scenarios.
The claim is proved.

 It follows that   $(u_1,v_1)=(x_1,y_1)\in\mathcal{P}'_{rs}$. Therefore, the first $k-1$ states  of the itinerary of $(\bu,\bv)$ are $id\to (S_{rs})^{k-2}$,  so $T_{M'}(\bu,\bv)\ge k-2=T_M(\bx,\by)-2$.

\medskip
\textit{Case 2. }$x_1s^{k-1}=\tau\gamma^{k-1}$ or $y_1r^{k-1}=\tau\gamma^{k-1}$.

Without loss of generality, we assume $x_1s^{k-1}=\tau\gamma^{k-1}$. Then $\tau\lhd_{r\gamma}y_1$ in $M$. By Lemma \ref{tool1} we have
$\bu=\kappa\alpha^{k-3}\cdots$ and $\bv=y_1r^{k-3}\cdots.$

If $r=\alpha$, then $\tau\lhd_{\alpha\gamma}y_1$, which forces $y_1=\kappa$ and $\bv=\kappa\alpha^{k-3}\cdots$.  Hence
 $T_{M'}(\bu,\bv)\geq |\bu\wedge\bv|\ge k-2$.  Similarly, if $r=\beta$, then $\tau\lhd_{\beta\gamma}y_1$, which forces $y_1=\lambda$ and $\bv=\lambda\beta^{k-3}\cdots$. So $T_{M'}(\bu,\bv)=T_{M'}(\kappa\alpha^{k-3}\cdots ,\lambda\beta^{k-3}\cdots)\ge k-2$.

If $M'$ is a $(\kappa,\tau)$ simplification of $M$,  we can prove  the lemma   in the same manner as above. The proof is finished.
\end{proof}

For $\bu,\bv\in\Omega$, we claim that if $U_1\ne V_1$, then
\begin{equation}\label{eq-key-2}
T_{M'}(\bu,\bv)-T_M(\bx,\by)\leq 5.
\end{equation}
Let $S_{M',1}$ be the first state of the itinerary of $(\bu,\bv)$ in $M'$ (after the initial state).
We will prove \eqref{eq-key-2}   in  Lemmas \ref{lem-key-5} and \ref{lem-key-6}.

\begin{lem}\label{tool2}
Let $\bv=v_1\eta^\ell\widetilde{\eta}\cdots$ where  $\widetilde{\eta}\neq \eta$.
If  $V_1=v_1$, then $g^{-1}(\bv)=v_1\eta^{\ell-2}\cdots$.
\end{lem}

\begin{proof}
Notice that if $\eta^p$ is a prefix of an element of ${\mathcal C}_{M'}$, then we must have  $p=1$
or $2$ (in the later case the element must be $\kappa\kappa\gamma$).
Hence,     we have
$|V_i|=1$ holds for $1\le i\le \ell-1$, so $g^{-1}(\bv)=v_1\eta^{\ell-2}\cdots$. 
\end{proof}

\begin{lem}\label{lem-key-5}
Equation \eqref{eq-key-2} holds if $U_1\neq V_1$ and $S_{M',1}=\text{Id}$.
\end{lem}

\begin{proof}
That $S_{M',1}=Id$ implies $u_1=v_1$, so at least one of $U_1$ and $V_1$ is in $\mathcal{C}_{M'}$, say $U_1\in\mathcal{C}_{M'}$. Now we divide the proof into two cases.

\medskip
\textit{Case 1.} $U_1=\tau\gamma\gamma$.

In this case, we have $|\bu\wedge\bv|\le 2$, so $T_{M'}(\bu,\bv)\leq 4$ and \eqref{eq-key-2} follows.

\medskip
\textit{Case 2.} $U_1=\kappa\alpha^k\kappa\gamma$ or $\kappa\alpha^k\kappa\gamma\gamma$($k\ge 0$).

\ding{172} If $V_1=\kappa\alpha^\ell\kappa\gamma$ or $\kappa\alpha^\ell \kappa\gamma\gamma$($\ell\ge 0$), then $\ell \neq k$ when $U_1$ and $V_1$ have the same form.  It is easy to see that
$
|\bu\wedge\bv|\leq 3+\min\{k,\ell\}.
$

Applying $g_0^{-1}$ to $U_1$ and $V_1$, we have
$$
X_1\in\{\tau\gamma^{k+2},\kappa\alpha^{k+1}\kappa\gamma\} \text{ and }
Y_1\in\{\tau\gamma^{\ell+2},\kappa\alpha^{\ell+1}\kappa\gamma\}.
$$
Therefore, $T_M(\bx,\by)\ge 1+\min\{k+1,\ell+1\}$. Thus
$$T_{M'}(\bu,\bv)-T_M(\bx,\by)\leq |\bu\wedge\bv|+2-T_M(\bx,\by)\leq 3.$$


\ding{173} If $V_1=\kappa$, write $\bv$ as $\kappa\alpha^\ell\widetilde{\alpha}\dots$ where $\widetilde{\alpha} \neq \alpha$, then
$
|\bu\wedge\bv|\leq 2+\min\{k,\ell\}.
$

 Applying $g_0^{-1}$ to $U_1$ and using  Lemma \ref{tool2} to $\bv$, we get
$$
X_1\in\{\tau\gamma^{k+2},\kappa\alpha^{k+1}\kappa\gamma\},\quad \by=\kappa\alpha^{\ell-2}\cdots.
$$
Therefore, $T_M(\bx,\by)\ge 1+\min\{k+1,\ell-2\}$. So \eqref{eq-key-2} holds by the same reason as above. The lemma is proved.
\end{proof}

\begin{lem}\label{lem-key-6}
Equation \eqref{eq-key-2} holds if $U_1\neq V_1$ and $S_{M',1}\in Q\setminus\{Id,Exit\}$.
\end{lem}

\begin{proof}
Let $k=T_{M'}(\bu,\bv)$. If $k\le 5$,  \eqref{eq-key-2} holds trivially. So we assume $k\ge 6$. Denote $S_{M',1}=S_{rs}$, where $r,s\in\{\alpha,\beta,\gamma\}$. Then the itinerary of $(\bu,\bv)$ is $id\to (S_{rs})^{k}\to Exit$, so $(u_1,v_1)\in\mathcal{P}_{rs}$ and
$$
\bu=u_1s^{k-1}\dots,\ \bv=v_1r^{k-1}\dots.
$$
Since $\tau$ is double-maximal in $M'$ (Lemma \ref{max_element}), we have $u_1s^{k-1}, v_1r^{k-1}\ne\tau\gamma^{k-1}$.

\medskip
\textit{Case 1. }$u_1s^{k-1}, v_1r^{k-1}\ne\kappa\alpha^{k-1}$.

In this case, neither  $\bu$ nor $\bv$ can be initialled by $\tau\gamma$ or $\kappa\alpha$, so $|U_1|=|V_1|=1$.
 Hence, by Lemma \ref{tool2} we have
$\bx=u_1s^{k-3}\cdots,\ \by=v_1r^{k-3}\cdots.$
Thus
$$T_M(\bx,\by)\geq T_{M'}(\bx,\by)\ge k-2=T_{M'}(\bu,\bv)-2.$$

\medskip
\textit{Case 2. }$u_1s^{k-1}=\kappa\alpha^{k-1}$ or $v_1r^{k-1}=\kappa\alpha^{k-1}$.

Without loss of generality, assume that $u_1s^{k-1}=\kappa\alpha^{k-1}$, then
$v_1\lhd_{\alpha r}\kappa$.
 Since $\kappa$ is $\alpha\gamma$-isolated in $M'$ (Lemma \ref{max_element}), we have $v_1=\lambda$ and $r=\beta$. So $M'$ is a $(\tau,\kappa, \lambda)$-simplification
of $M$.
On one hand,
$\by=g^{-1}(\lambda\beta^{k-1}\cdots)=\lambda\beta^{k-1}\cdots$.
On the other hand, if $U_1=\kappa\alpha^\ell\kappa\gamma(\ell\ge k-1)$, then
$X_1=g_0^{-1}(U_1)=\tau\gamma^{\ell+2}$; otherwise, by Lemma \ref{tool2} we have $\bx=\kappa\alpha^{k-3}\cdots$
no matter $U_1=\kappa\alpha^\ell\kappa\gamma\gamma(\ell\ge k-1)$ or $|U_1|=1$. Thus
$$
\bx\in\{\tau\gamma^{k+1}\cdots,\kappa\alpha^{k-3}\cdots\} \text{ and } \by=\lambda\beta^{k-1}\dots,
$$
which imply that $T_M(\bx,\by)\ge k-2$. The lemma is proved.
\end{proof}

\begin{proof}[\textbf{Proof of Theorem \ref{main}}]
For  $\bx,\by\in\Omega$, if $X_1\dots X_k=Y_1\dots Y_k$ and $X_{k+1}\ne Y_{k+1}$, then $U_1\dots U_k=V_1\dots V_k$ and $U_{k+1}\ne V_{k+1}$. Let $\ell=|X_1\dots X_k|$, then
\begin{equation}\label{proofofmain1}
T_M(\bx,\by)-T_{M'}(\bu,\bv)=T_M(\sigma^\ell(\bx),\sigma^\ell(\by))-
T_{M'}(\sigma^\ell(\bu),\sigma^\ell(\bv)),
\end{equation}
where $\sigma$ is the shift operator. By \eqref{eq-key-1} and \eqref{eq-key-2},  Theorem \ref{main} holds for $\sigma^\ell(\bx)$ and $\sigma^\ell(\by)$,  so it also holds for $\bx$ and $\by$.
\end{proof}

\medskip

\begin{appendix}

\section{\textbf{Connected components of a fractal gasket}}\label{app:A}

Define $\pi(x,y)=y$.

\begin{lem}\label{component}
Let $K$ be a fractal gasket. If $K$ satisfies the top isolated condition or $\omega_\gamma\not\in K$, then any nontrivial connected component of $K$ is a line segment. If $\alpha$ and $\beta$ are not in the same horizontal block in addition, then $K$ is totally disconnected.
\end{lem}

\begin{proof}
Let $U$ be a   connected component of $K$. Notice that the collection of the vertices of $\varphi_i(\triangle)$, $i\in\Sigma$, forms a cut set of $K$.

Let $I$ be a horizontal block of $\Sigma$ (see Definition \ref{Hblock}), we denote
$B(I)=\bigcup_{i\in I} \varphi_i(\triangle)$.

If $\omega_\gamma\not \in K$, then all $B(I)$ are disjoint, so either $U\subset B(I)$ or they are disjoint;
if $K$ satisfies the top isolated condition, it is easy to show that the same conclusion still holds.
  It follows that $|\pi(U)|\leq r^*$.

   Now we regard $K$ as the invariant set of the IFS
$\{\varphi_I\}_{I\in\Sigma^k}$. By the same argument as above, we obtain $|\pi(U)|\leq (r^*)^k$.
Let $k \rightarrow\infty$, we see that $\pi(U)$ is a single point. So $U$ is a horizontal line segment.

 Suppose $\alpha$ and $\beta$ are not in the same horizontal block, let $p$ be the size of the block containing $\alpha$, and $q$ be the size of the block containing $\beta$. Then for any $k\geq 1$, the size of any horizontal block
 of $\Sigma^k$ is bounded by $p+q$. Hence the diameter of $U$ is bounded by $(p+q)(r^*)^k$. Let $k\to \infty$,
 we obtain that $U$ is a single point.
\end{proof}

\section{The transducer of $g$}
\begin{figure}
  \centering
  \includegraphics[width=12 cm]{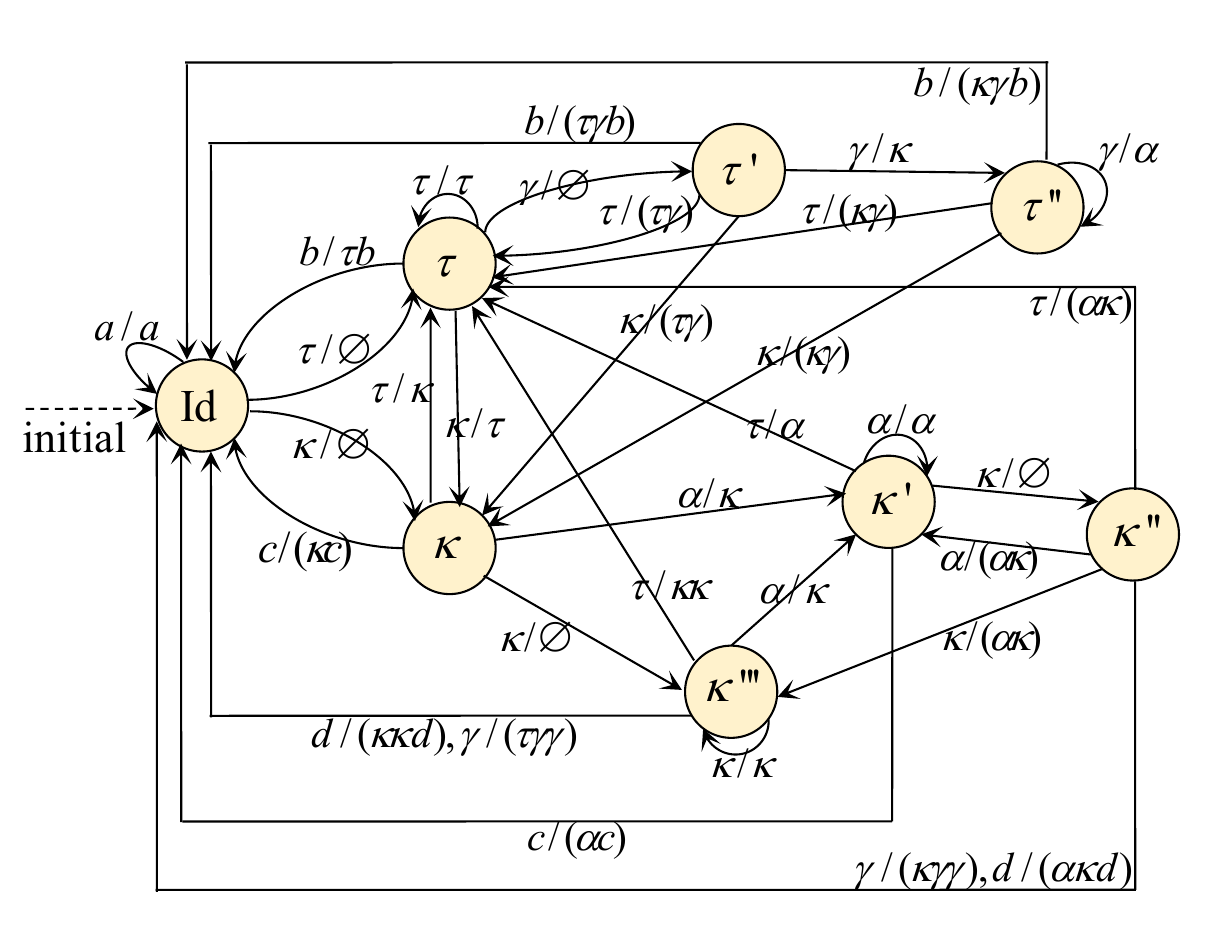}\\
\caption{The transducer of $g$, where $a\in\Sigma\setminus\{\tau,\kappa\}$, $b\in\Sigma\setminus\{\tau,\kappa,\gamma\}$, $c\in\Sigma\setminus\{\tau,\kappa,\alpha\}$, $d\in\Sigma\setminus\{\tau,\kappa,\alpha,\gamma\}$.}
\label{transducer}
\end{figure}

The map $g$ defined in \eqref{gxdecomposition} can be realized by the transducer indicated in Figure \ref{transducer}.
The state set is $\{Id,\tau,\tau',\tau'',\kappa,\kappa',\kappa'',\kappa'''\}$ where `$Id$' is the initial state. Each edge is labeled by $x_i/\omega_i$, where $x_i\in\Sigma$ is the input letter, and $\omega_i\in\Sigma^*$ is the output word. For any input symbol string $\bx=(x_i)_{i=1}^\infty\in\Sigma^\infty$, there is a unique sequence determined by the transducer, which is denoted by $\omega_1\omega_2\dots=g(\bx)$.

\end{appendix}


\begin{thebibliography}{99}
\addcontentsline{toc}{chapter}{Bibliography}

\bibitem {Bal2000}
R. Balakrishnan and K. Ranganathan, A Textbook of Graph Theory, \emph{Springer-Verlag, Berlin Heidelberg,} 2000.

\bibitem{BanM09}
C. Bandt and M. Mesing: Self-affine fractals of finite type, \textit{Banach Center Publication}, \textbf{84} (2009) 131-148.

\bibitem{Bao18}
T.Q. Bao, S. Cobza, A. Soubeyran: Variational principles, completeness and the existence of traps in behavioral sciences, \textit{Annals of Operations Research}, \textbf{269} (2018), 53-79.

\bibitem{Bonk} M. Bonk and S. Merenkov: Quasisymmetric rigidity of square Sierpinski carpets, \emph{Anal. Math.}, \textbf{177} (2013), no. 2, 591-643.

\bibitem{DS}
G. David and S. Semmes: Fractured fractals and broken dreams: self-similar geometry through metric and measure, \textit{Oxford Univ Press}, 1997.


\bibitem{FanRZ15}
A.H. Fan, H. Rao and Y. Zhang: Higher dimensional Frobenius problem: Maximal saturated cone, growth function and rigidity, \emph{J. Math. Pures Appl}, \textbf{104} (2015), 533-560.

\bibitem{JEH79}
J.E. Hopcroft, R. Motwani, Rotwani and J. D. Ullman: Introduction to Automata Theory, Language, and Computation, Addison-Wesley, Massachusetts (1979).

\bibitem{Hutchinson1981}
J.E. Hutchinson, Fractals and self-similarity, Indiana Univ. Math. J. 30 (1981) 713-747.

\bibitem{LuoL13}
J.J. Luo and K.S. Lau: Lipschitz equivalence of self-similar sets and hyperbolic boundaries, \textit{Adv. Math}, \textbf{235} (2013), 555-579.

\bibitem{LuoL16}
J.J. Luo and J.C. Liu: On the classification of frctal squares, \textit{Fractals}, \textbf{24} (2016), no. 1, 11 pp.

\bibitem{FM}
D. T. Marsh and K.J. Falconer: On the Lipschitz equivalence of Cantor sets,\textit{ Mathematika}, \textbf{39} (1992), 223-233.

\bibitem{Pepo90}
B. Pepo: Fixed points for contractive mappings of third order in pseudo-quasimetric spaces, \textit{Indag. Mathem. N. S.}, \textbf{4} (1990), no. 1,  473-482.

\bibitem{RWW17}
F. Rao, X.H. Wang , S.Y. Wen: On the topological classification of fractal squares, \textit{ Fractals}, \textbf{25} (2017), no. 3, 12 pp.

\bibitem{RRX06}
H. Rao, H.J. Ruan and L.F. Xi: Lipschitz equivalence of self-similar sets, \textit{C. R. Acad. Sci. Paris, Ser.I}, \textbf{342} (2006), 191-196.

\bibitem{RaoZ15}
H. Rao and Y. Zhang: Higher dimensional Frobenius problem and Lipschitz equivalence of Cantor sets. \textit{J. Math. Pures Appl}, \textbf{104} (2015), 868-881.

\bibitem{RaoZhu16}
H. Rao and Y.J. Zhu, \textit{Lipschitz equivalence of fractal squares and finite state automation}, arXiv preprint arXiv:1609.04271 (2016).

\bibitem{RuanW17} H. J. Ruan and Y. Wang: Topological invariants and Lipschitz equivalence of fractal sets,\textit{J. Math. Anal. Appl.}, \textbf{451} (2017), 327-344.

\bibitem{RuanWX14}
H. J. Ruan, Y. Wang and L.F. Xi: Lipschitz equivalence of self-similar sets with touching structures, \textit{Nonlinearity},  \textbf{27} (2014), no.6, 1299-1321.


\bibitem{Solomyak10}
K. I. Ero\v{g}lu, S. Rohde, B . Solomyak: Quasisymmetric conjugacy between quadratic dynamics and iterated function systems, \textit{Ergodic Theory $\&$ Dynamical Systems},  \textbf{30} (2010), no. 6, 1665-1684.

\bibitem{Why58}
G. T. Whyburn: Topological characterization of the Sierpinski curve, \textit{Fund. Math}, \textbf{45} (1958), 320-324.

\bibitem{XiXi10}
L.F. Xi and Y. Xiong: Self-similar sets with initial cubic patterns, \textit{C. R. Acad. Sci. Paris, Ser.I}, \textbf{348} (2010), 15-20.


\bibitem{XiXi20}
L.F. Xi and  Y. Xiong: Algebraic criteria for Lipschitz equivalence of dust-like self-similar sets, \textit{J.  London Math. Soc.}, (2020)
(DOI: 10.1112/jlms.12392).

\bibitem{YZ18}
Y.J. Zhu and Y.M. Yang: Lipschitz equivalence of self-similar sets with two-state automation, \textit{J. Math. Anal. Appl}, \textbf{458} (2018), no. 1, 379-392.











\end{thebibliography}
\end{document}